\documentclass[10pt,twoside,english]{article}
\usepackage{amssymb,amsmath,babel,geometry,latexsym,graphics,tabularx,shapepar,enumerate}
\usepackage{hyperref}
\usepackage[all,2cell]{xy} \UseAllTwocells \SilentMatrices
\def\carlitz{{{\text{Car}}}}

\newtheorem{Theorem}{Theorem}
\newtheorem{Lemme}[Theorem]{Lemma}
\newtheorem{Proposition}[Theorem]{Proposition}

\newtheorem{Definition}[Theorem]{Definition}
\newtheorem{Remarque}[Theorem]{Remark}

\newcommand{\ol}{\overline}
\newcommand{\A}{\mathbb{A}}
\newcommand{\ZZ}{\mathbb{Z}}
\newcommand{\FF}{\mathbb{F}}
\newcommand{\CC}{\mathbb{C}}

\newcommand{\GG}{\mathbb{G}}

\newcommand{\LL}{\mathbb{L}}
\newcommand{\MM}{\mathbb{M}}

\newcommand{\TT}{\mathbb{T}}

\newcommand{\hPsi}{\boldsymbol{\widehat{\Psi}}}
\newcommand{\tg}{\widetilde{g}}
\newcommand{\tDelta}{\widetilde{\Delta}}

\newcommand{\bfs}{\boldsymbol{s}}
\newcommand{\unif}{u}

\newcommand{\bsb}{\boldsymbol}

\newcommand{\sqm}[4]{
\left(\begin{array}{ll}#1 & #2 \\ #3 & #4\end{array}\right)}
 
\newcommand\CVD{{\hfill\hfil{\lower 2 pt\hbox{\vrule\vbox to 7pt 
{\hrule width 6pt\vfill\hrule}\vrule}}}\vskip 0.5cm}

\let\le=\leqslant  
\let\ge=\geqslant                                 

\title{Estimating the order of vanishing at infinity \\ of Drinfeld quasi-modular forms.}
\author{Federico Pellarin}

\begin{document}

\maketitle

\begin{small}
\noindent\textbf{Abstract.}
We introduce and study certain deformations of 
{\em Drinfeld quasi-modular forms} by using {\em rigid analytic trivialisations} of corresponding {\em Anderson's $t$-motives}.
We show that a sub-algebra of these deformations has a natural graduation by the group $\ZZ^2\times\ZZ/(q-1)\ZZ$ and 
an homogeneous automorphism,
and we deduce from this and other properties {\em multiplicity estimates}. 

%
\end{small}

\medskip

\setcounter{tocdepth}{1}
\tableofcontents

\section{Introduction, results.}

Let $E_4,E_6$ be the Eisenstein series of weights $4,6$ respectively, rescaled so that
they have limit $1$ as the imaginary part of the variable tends to infinity; let us denote by $M$ the two-dimensional $\CC$-algebra $\CC[E_4,E_6]$.
We have $M=\oplus_wM_w$, where $M_w$ denotes 
the $\CC$-vector space of classical modular forms $\mathbf{SL}_2(\ZZ)$ of weight $w$. 

Let $E_2$ be the so-called ``false" Eisenstein series of weight $2$ (rescaled), and
let us consider, for $l,w$ non-negative integers with $w>0$ even,
the non-zero finite dimensional $\CC$-vector space
$$\widetilde{M}_{w}^{\leq l}=M_w\oplus M_{w-2}E_2\oplus\cdots\oplus E_2^lM_{w-2l}$$
of {\em classical quasi-modular forms} of weight $w$ and depth $\leq l$.

The local behaviour at infinity of quasi-modular forms yields
an embedding $\widetilde{M}_{w}^{\leq l}\rightarrow\CC[[q]]$ where $q=e^{2\pi iz}$, $z\in\mathcal{H}$ being the variable in the complex upper-half plane.
Let $f$ be in $\widetilde{M}_{w}^{\leq l}\setminus\{0\}$. Then,
we may write $f=q^{\nu_\infty(f)} g$ with $g$ a unit of $\CC[[q]]$, for an integer $\nu_\infty(f)\geq 0$ which is uniquely determined; this is the {\em order of vanishing at infinity}
of $f$. 

A simple resultant argument suffices to show that
\begin{equation}\label{boundwith dimension}
\nu_\infty(f)\leq 3\dim_\CC\widetilde{M}^{\leq l}_w.
\end{equation}

Now, let $q=p^e$ (\footnote{The double use of the letter $q$ in this paper will not be a source of confusion; other harmless abuses of notation will appear in this paper.}) be a power of a
prime number $p$ with $e>0$ an integer, let $\FF_q$ be the finite field with $q$ elements.
Let us write $A=\FF_q[\theta]$ and $K=\FF_q(\theta)$, with $\theta$
an indeterminate over $\FF_q$, and define an absolute value
$|\cdot|$ on $K$ by $|a|=q^{\deg_\theta a}$, $a$ being in $K$, so that
$|\theta| = q$.  Let $K_\infty :=\FF_q((1/\theta))$ be the
completion of $K$ for this absolute value, let $K_\infty^{\text{alg.}}$ be an
algebraic closure of $K_\infty$, let $C$ be the completion of
$K_\infty^{\text{alg.}}$ for the unique extension of $|\cdot|$ to $K_\infty^{\text{alg.}}$, and let $K^{\text{alg.}}$ be the algebraic closure of $K$
in $C$. 

Following Gekeler in \cite{Ge}, we denote by $\Omega$ the rigid analytic space $C\setminus K_\infty$ and write $\Gamma$ for $\mathbf{GL}_2(A)$, group that acts on $\Omega$ by homographies. In this setting we have three functions $E,g,h:\Omega\rightarrow C$, holomorphic in the sense of \cite[Definition 2.2.1]{FP},
such that, for all $\gamma=\begin{pmatrix}a&b\\ c&d\end{pmatrix}\in\Gamma$ and $z\in\Omega$:
\begin{eqnarray}
g(\gamma(z))&=&(cz+d)^{q-1}g(z),\nonumber\\
\quad h(\gamma(z))&=&(cz+d)^{q+1}\det(\gamma)^{-1}h(z),\nonumber\\
E(\gamma(z))&=&(cz+d)^2\det(\gamma)^{-1}\left(E(z)-\frac{c}{\widetilde 
{\pi}(cz+d)}\right)\label{formE}
\end{eqnarray}
where $\gamma(z)=(az+b)/(cz+d)$ and $\widetilde{\pi}$ is a fundamental period of {\em Carlitz module}, defined by the convergent product:
$$\widetilde{\pi}:=\theta(-\theta)^{\frac{1}{q-1}}\prod_{i=1}^\infty(1-\theta^{1-q^i})^{-1}\in K_\infty((-\theta)^{\frac{1}{q-1}})\setminus K_\infty,$$ a choice of a $(q-1)$-th root of $-\theta$ having been made once and for all (\footnote{See \cite[Section 2.1]{Bourbaki}, where the notation $\overline{\pi}$ is adopted; there is an analogy with the number $2\pi i$.}).

The functional equations above tell that $g,h$ are {\em Drinfeld modular forms},
of {\em weights} $q-1$, $q+1$ and {\em types} $0,1$ respectively. For $w$ integer and $m\in\ZZ/(q-1)\ZZ$, we denote by $M_{w,m}$ the $C$-vector space 
of Drinfeld modular forms of weight $w$ and type $m$.

After (\ref{formE}), the function $E$ is not a Drinfeld modular form. In \cite{Ge},
Gekeler calls it ``False Eisenstein series" of weight $2$ and type $1$; it is easy to show that $E,g,h$ are algebraically independent.
For $l,w$ non-negative integers and $m$ a class of $\ZZ/(q-1)\ZZ$, we introduce the $C$-vector space of {\em Drinfeld quasi-modular forms
of weight $w$, type $m$ and depth $\leq l$}:
$$\widetilde{M}^{\leq l}_{w,m}=M_{w,m}\oplus M_{w-2,m-1}E\oplus\cdots\oplus M_{w-2l,m-l}E^l.$$

Let $e_\carlitz:C\rightarrow C$ be the {\em Carlitz exponential function} (see below, (\ref{exponential}))
and let us write $u:\Omega\rightarrow C$ for the ``parameter at infinity'' of $\Omega$, that is, the function
$u(z)=1/e_\carlitz(\widetilde{\pi}z)$; the $C$-algebra $\widetilde{M}$ embeds in $C[[u]]$ (cf. \cite{Ge}). 

If $w,m,l$ are such that the finite dimensional vector space
$\widetilde{M}^{\leq l}_{w,l}$ does not vanish, any $f$ non-zero in $\widetilde{M}^{\leq l}_{w,l}$ can be written, for $u=u(z)$, as
$f(z)=u^{\nu_\infty(f)}\phi$ with $\phi$ a unit of the ring $C[[u]]$ (convergence when $|u|$ is small enough) for some non-negative integer
$\nu_\infty(f)$ uniquely determined; the order of vanishing at infinity of $f$. 

The aim of this paper is to prove the following analog of (\ref{boundwith dimension}):
\begin{Theorem}\label{secondtheorem} 
Let $l$ be a positive integer.
There exist two constants $c_1,c_2$, with $c_1$ depending on $l,q$ and $c_2$ depending on $q$, with the following property.
If $w\geq c_1$ and $\widetilde{M}^{\leq l}_{w,m}\neq(0)$ then, for all $f\in \widetilde{M}^{\leq l}_{w,m}\setminus\{0\}$,
$$\nu_\infty(f)\leq c_2\dim_C\widetilde{M}^{\leq l}_{w,m}.$$
\end{Theorem}

Our result is completely explicit; if $l>0$ and
$f\in\widetilde{M}^{\leq l}_{w,m}$ is a non-zero Drinfeld quasi-modular form, we will prove that
\begin{equation}\label{finalupperbound}\nu_\infty(f)\leq16q^3(3 + 2 q)^2lw,
\end{equation}
provided that
\begin{equation}\label{condition}w\geq 4l\left(2 q(q+2)(3 + 2 q)l+3(q^2+1)\right)^{3/2},\end{equation}
which obviously implies the Theorem by a simple computation of the dimension of $\widetilde{M}^{\leq l}_{w,m}$.

In \cite{BP} we conjectured that for all $z\in\Omega$ at least three of the four numbers
$u(z),E(z),g(z),h(z)$
are algebraically independent over $K$. If true, this statement would be an analog of Nesterenko's celebrated theorem 
on the algebraic independence of values of Eisenstein's series and the parameter at infinity 
(see \cite[Chapter 3, Theorem 1.1]{NP} and \cite{Nes}). In the attempt of proving this conjecture,
it turned out very difficult to adopt Nesterenko's original scheme of proof. Indeed, Drinfeld quasi-modular 
forms, while sharing several superficial similarities with classical quasi-modular forms, essentially behave in a
different way. The main difficulties encountered are the following:
\begin{description}
\item Absolute values of cofficients of $u$-expansions of Drinfeld quasi-modular forms grow ``too rapidly" depending on the index;
if $f=\sum_ic_iu^i$ is such a form, then the estimate $|c_i|\ll e^{ci}$ is best possible (in the classical case,
we would have $|c_i|\ll i^c$).
\item The algebra $\widetilde{M}$ is endowed with higher derivatives (this was studied in the joint work \cite{BP}); however,
this structure alone does not seem to easily deliver a suitable analogue of the {\em separation property} 
\cite[Lemma, p. 212]{BM}, useful for multiplicity estimates in differential algebras. 
\end{description}
These difficulties suggest that the algebra $\widetilde{M}$ {\em is not} an appropriate environment to 
study the arithmetic of values of Drinfeld quasi-modular forms. 
Theorem \ref{secondtheorem} could superficially look like a mere copy of the elementary inequality 
(\ref{boundwith dimension}). Before our proof, it was however very resistant to any attempt to prove it.
The main motivation of this paper is to find new structures allowing to prove Theorem \ref{secondtheorem}.

In this paper, we introduce a new class of functions (deformations of Drinfeld quasi-modular forms), endowed 
with certain automorphic properties and a ``Frobenius structure". 
The theory we introduce is strongly influenced by that of {\em Anderson's $t$-motives}. The idea of appealing to 
Anderson's theory is very natural; however, a new entity occurs here, making this paper useful:  the functions we deal with
have automorphic properties; they generate an algebra graded by the group $\ZZ^2\times\ZZ/(q-1)\ZZ$. Similar 
objects in the classical theory do not seem to be already known.

The main properties of our functions
are listed in Proposition \ref{firsttheorem} below for the sake of commodity (this proposition will not be applied directly).
With the help of {\em all} these properties and a {\em transcendence proof}, we deduce Theorem \ref{secondtheorem}.


\medskip

In order to present Proposition \ref{firsttheorem} we require some further preparation. 

Let $t$ be an independent indeterminate, let us temporarily denote by $\TT$ the subring of formal series of $C[[t]]$ converging for all $t\in C$ such that $|t|<r$ for a certain real number $r>q$
(\footnote{Later, we will see that $r\geq q^q$ and we will then use the notation $\TT_{<q^q}$.}).
In Section \ref{bsbE} we will construct a function $$\bsb{E}:\Omega\rightarrow\TT$$
such that, for $t\in C$ and $z\in\Omega$ with $|t|$ and $|u|$ small enough (and with $u=u(z)$),
the value $\bsb{E}(z)(t)$ is that of a convergent double series in $C[[t,u]]$. More precisely,
we will show (Proposition \ref{propositionEuexpansion}) the existence of polynomials $c_i\in\FF_q[t,\theta]$ such that there is a locally convergent expansion
\begin{equation}
\bsb{E}(z)(t)=u\sum_{i\geq 0}c_i(t)u^{i(q-1)}\in u\FF_q[t,\theta][[u^{q-1}]].
\label{bsbEidentifieswithformalseries}\end{equation} We identify $\bsb{E}$ with the formal series at the right-hand side of (\ref{bsbEidentifieswithformalseries}).

We will use the following extension of {\em Anderson's} $\FF_q[t]$-linear map $\tau:C[[t,u]]\rightarrow C[[t,u]]$:
$$\tau\sum_{m,n\geq 0}c_{m,n}t^mu^n:=\sum_{m,n\geq 0}c_{m,n}^qt^mu^{qn}$$ the $c_{m,n}$'s being in $C$.

Let $\bsb{F}$ be the formal series $\tau\bsb{E}\in u\FF_q[t,\theta][[u^{q(q-1)}]]$; it converges locally at $(0,0)$, and extends to a well defined function
$\Omega\rightarrow\TT$. Let us denote by  $\MM^\dag$ the algebra
$\TT[g,h,\bsb{E},\bsb{F}]$, that we will often identify with its image in $C[[t,u]]$. We have the following Proposition (\footnote{
It results from the combination of the six Propositions
\ref{taudifferenceeq}, \ref{propositionautomorphyE}, \ref{propositionEuexpansion}, \ref{gradingfiltering}, \ref{algebraicindependenceEghE1} and \ref{almost_bijective},
and elementary considerations. The Proposition is stated to ease the access of the paper, but later, we will require the full statement of the six propositions above. Throughout these six propositions, many other properties of the functions $\bsb{E},\bsb{F}$ and of the 
algebra $\MM^\dag$ will be highlighted.}):

\begin{Proposition}\label{firsttheorem}
The algebra $\MM^\dag\subset C[[t,u]]$ enjoys the following properties.
\begin{enumerate}
\item The dimension of $\MM^\dag$ is four, so that the formal series $g,h,\bsb{E},\bsb{F}$ are algebraically independent over $C((t))$ and define a basis
of $\MM^\dag$.
\item The basis $(g,h,\bsb{E},\bsb{F})$ of $\MM^\dag$ above is constituted by formal series in $\FF_q[t,\theta][[u]]$ (it is integral over $\FF_q[t,\theta]$).
\item The algebra $\MM^\dag$ is graded by the group $G=\ZZ^2\times\ZZ/(q-1)\ZZ$. In other words, $\MM^\dag=\oplus_{(\mu,\nu,m)\in G}\MM^\dag_{\mu,\nu,m}$.
We further have $\MM^\dag_{0,0,0}=\TT$ and $\MM^\dag_{\mu,0,m}=M_{\mu,\nu}\otimes_C\TT$. For this graduation,
the ``degrees" of $g,h,\bsb{E}$ and $\bsb{F}$ are respectively the following elements of $G$: $(q-1,0,0),(q+1,0,1),(1,1,1)$ and $(q,1,1)$.
\item For any $\bsb{f}\in\MM^\dag_{\mu,\nu,m}$, the formal series $\varepsilon(\bsb{f})=\bsb{f}|_{t=\theta}$ is a well defined Drinfeld 
quasi-modular form of $\widetilde{M}^{\leq \nu}_{\mu+\nu,m}$ and we have a surjective $C$-linear map $\varepsilon:\MM^\dag\rightarrow\widetilde{M}$
such that $\varepsilon(\bsb{E})=E$; in this sense, $\bsb{E}$ is a deformation of $E$.
\item The algebra $\MM^\dag$ is stable under the action of $\tau$: more precisely, $\tau$ induces $\FF_q[t]$-linear maps
$\MM^\dag_{\mu,\nu,m}\rightarrow\MM^\dag_{q\mu,\nu,m}$.
\item Let $\bsb{f}$ be in $\MM^\dag$, let us assume that $\bsb{f}=\sum_{i\geq 0}c_iu^i$ with $c_i\in\FF_q[t,\theta]$.
Then, $\deg_tc_i\leq c\log i$ where $c$ is a constant depending on $\bsb{f}$
only.
\end{enumerate}
\end{Proposition}

The properties described in Proposition \ref{firsttheorem}, a variant of Siegel's Lemma (Proposition \ref{functionsfmunum})
and a transcendence construction will be used to prove Theorem \ref{secondtheorem}.

\section{Anderson's functions\label{section:trivialisations}}
 
In this section we recall some tools introduced in \cite[Section 2]{Ge}, \cite{An, ABP} and described in  \cite[Section 2 and Section 4.2]{Bourbaki}.
 
As {\em $A$-lattice of rank $r>0$} we mean a free sub-$A$-module of $C$ of rank $r$, discrete in the sense that, given a compact subset of $C$, only finitely many elements of it lie in.
Let $\Lambda\subset C$ be an $A$-lattice of rank $r$
and let us consider, for $\zeta\in C$, the {\em exponential function associated to $\Lambda$}, defined by the product: 
\begin{equation}\label{weierstrassproduct}
e_\Lambda(\zeta):=\zeta\prod_{\omega\in\Lambda\setminus\{0\}}\left(1-\frac{\zeta}{\omega}\right),
\end{equation}
which converges for all $\zeta\in C$.
For $\lambda\in C^\times$, the product expansion (\ref{weierstrassproduct}) implies:
\begin{equation}\label{homogeneity}
e_{\lambda\Lambda}(\zeta)=\lambda e_{\Lambda}(\lambda^{-1} \zeta).
\end{equation}
There exist elements $1=\alpha_0(\Lambda),\alpha_1(\Lambda),\alpha_2(\Lambda),\ldots\in C$, depending on $\Lambda$ only, such that:
\begin{equation}\label{exponentialfunction}
e_\Lambda(\zeta)=\sum_{n\geq 0}\alpha_n(\Lambda)\zeta^{q^n},
\end{equation}
the series having infinite radius of convergence (cf. \cite{Ge, Go2}).

The construction of the exponential function by (\ref{weierstrassproduct}) is the main tool to prove that the category whose objects are homothecy classes of $A$-lattices of rank $r$ and morphisms are inclusions, is dually equivalent to the category whose objects are isomorphism classes of Drinfeld $A$-modules of rank $r$ and morphisms are isogenies
(see \cite[Section (2.6)]{Ge} or \cite[Section 2]{Bourbaki}). For $\Lambda$ as above, 
there is a Drinfeld $A$-module $\phi_\Lambda$ such that
\begin{equation}\label{drinfeldmodule}
\phi_\Lambda(a)e_\Lambda(\zeta)=e_\Lambda(a\zeta)
\end{equation}
(for all $\zeta\in C$ and $a\in A$), which is uniquely determined by its value $\phi_\Lambda(\theta)\in\mathbf{End}_{\FF_q-\text{lin.}}(\GG_a(C))$ in $\theta$. This value
is a polynomial of degree $r$ in $\tau$, which we recall, is the Frobenius endomorphism $\tau:c\mapsto c^q$.
On the other side, to any Drinfeld $A$-module $\phi$ of rank $r$, a lattice $\Lambda_\phi$ of rank $r$ can be associated, so that the functors $\Lambda\mapsto\phi_\Lambda$
and $\phi\mapsto\Lambda_\phi$ are inverse of each other up to isomorphisms.

Let $t$ be a new indeterminate. 
With $\Lambda$ an $A$-lattice of rank $r>0$ and $e_\Lambda$ as in (\ref{weierstrassproduct}), 
let us consider $\omega\in\Lambda\setminus\{0\}$ and
introduce, following Anderson in \cite{An}, the formal series:
$$s_{\Lambda,\omega}(t):=\sum_{i=0}^\infty e_\Lambda\left(\frac{\omega}{\theta^{i+1}}\right)t^i.$$

For a positive real number $r$, we denote by $\TT_{<r}$ the sub-$C$-algebra of $C[[t]]$ whose elements are formal series $\sum_{i\geq 0}c_it^i$ that converge 
for any $t\in C$ with $|t|<r$. We denote by $\TT_{>0}$ the sub-$C$-algebra of $C[[t]]$ whose series converge in some open disk containing $0$, and we notice
 that all the series of $\TT_{<q^q}$ converge at $t=\theta$. We also denote by $\TT_\infty$ the sub-$C$-algebra of series that converge everywhere in $C$.

If $r_1>r_2>0$, we have $$\TT_{>0}\supset\TT_{<r_2}\supset\TT_{< r_1}\supset\TT_\infty.$$ The {\em Tate algebra} 
of formal series of $C[[t]]$ converging for all $t$ such that $|t|\leq 1$ will be denoted by $\TT$; it is contained in $\TT_{<1}$ and contains
$\TT_{<1+\epsilon}$ for all $\epsilon>0$.

It is easy to verify that, with $\Lambda$ and $\omega\in\Lambda$ as above,
$s_{\Lambda,\omega}\in\TT_{<q}\subset\TT$. If $\Lambda\subset K_\infty^{\text{alg.}}$,
it can be proved that $s_{\Lambda,\omega}(t)\in K_\infty^{\text{alg.}}[[t]]$.

We extend the operator $\tau$ from $C$ to $C[[t]]$ as follows:
$$f=\sum_{n\geq 0}c_nt^n\mapsto \tau f:=\sum_{n\geq 0}c_n^qt^n.$$
We will also write $f^{(k)}$ for $\tau^kf$, $k\in\ZZ$ (the operator $\tau^{-1}$ is well defined). 
One checks that $\tau$ sends $\TT_{<r}$ in $\TT_{<r^q}$. The extension $\tau$ so constructed defines $\FF_q$-automorphisms of $\TT_{>0},\TT$ and $\TT_\infty$.

We write $\overline{A}=\FF_q[t],\overline{K}=\FF_q(t)$.
If $a=a(\theta)\in A$ we also write $\overline{a}=a(t)\in\overline{A}$.
If $\Lambda$ is an $A$-lattice
of rank $r$ and if $\phi_\Lambda$ is the Drinfeld $A$-module of rank $r$ in (\ref{drinfeldmodule}), then, for all $a_1,a_2\in A$
and $\omega_1,\omega_2\in\Lambda$,
\begin{equation}\label{pseudolinearity}
\phi_\Lambda(a_1)s_{\Lambda,\omega_1}+\phi_\Lambda(a_2)s_{\Lambda,\omega_2}=s_{\Lambda,a_1\omega_1+a_2\omega_2}=\overline{a}_1s_{\Lambda,\omega_1}+\overline{a}_2s_{\Lambda,\omega_2}.
\end{equation}
These identities, which hold in $\TT$, are proved in \cite[Section 4.2.2]{Bourbaki}.

From (\ref{homogeneity}) it immediately follows that, for $\lambda\in C^\times$,
\begin{equation}\label{homogeneity2}
s_{\lambda \Lambda,\lambda\omega}(t)=\lambda s_{\Lambda,\omega}(t).
\end{equation}
We also have the series expansion (cf. \cite[Section 4.2.2]{Bourbaki})
\begin{equation}\label{residues}
s_{\Lambda,\omega}(t)=\sum_{n=0}^\infty\frac{\alpha_n(\Lambda)\omega^{q^n}}{\theta^{q^n}-t},
\end{equation}
uniformly convergent in every compact subset of $C\setminus\{\theta,\theta^q,\ldots\}$, and
$s_{\Lambda,\omega}(t)-\omega/(\theta-t)$ extends to a rigid holomorphic function for $|t|<q^q$. We will then
often say that $s_{\Lambda,\omega}$ has a {\em simple pole of residue} $-\omega$ in $t=\theta$. Notice that other poles occur at $t=\theta^q,\theta^{q^2},\ldots$, but we will never need to focus on them in this paper.

\medskip\noindent\emph{Example: rank one case.}
For $\Lambda=\widetilde{\pi}A$ (rank $1$),
the exponential function (\ref{exponentialfunction}) is:
\begin{equation}\label{exponential}
e_\carlitz(\zeta)=\sum_{n\ge 0}\frac{\zeta^{q^n}}{d_n},
\end{equation}
where $d_0:=1$ and $d_i:=[i][i-1]^q\cdots[1]^{q^{i-1}}$, recalling that $[i]=\theta^{q^i}-\theta$ if $i>0$.
The relations (\ref{drinfeldmodule}) become, for all $a\in A$, $$\phi_\carlitz(a)e_\carlitz(\zeta)=e_\carlitz(a\zeta),$$
where $\phi_\carlitz$ is Carlitz's module defined by $$\phi_\carlitz(\theta)=\theta\tau^0+\tau\in\mathbf{End}_{\FF_q\text{-lin.}}(\GG_a)$$
(see Section 4 of \cite{Ge}).

We will write $s_\carlitz=s_{\widetilde{\pi}A,\widetilde{\pi}}$.
The function $s_\carlitz$ has a simple pole in $\theta$ with residue $-\widetilde{\pi}$.

By (\ref{pseudolinearity}) (cf. \cite[Section 4.2.5]{Bourbaki}), the following $\tau$-difference equation holds:
\begin{equation}\label{taudifferencebfs}
s_\carlitz^{(1)}(t)=(t-\theta)s_\carlitz.
\end{equation}

After \cite[Theorem 2.2.9]{FP}, $\TT$ is a principal ideal domain. This property can be used to verify that
the subfield of constants $\LL^\tau:=\{l\in\LL,\tau l=l\}$, where $\LL$ is the fraction field of $\TT$, is equal to $\ol{K}:=\FF_q(t)$ (see also \cite[Lemma 3.3.2]{Pa}). We deduce, just as in the proof of \cite[Lemma 3.3.5]{Pa}, that the $\tau$-difference equation $f^{(1)}=(t-\theta)f$ has, as a complete set of solutions in $\LL$, the $\FF_q(t)$-vector space $\FF_q(t)s_\carlitz$. In fact, for all $a=a(\theta)\in A$, we have $s_{\widetilde{\pi}A,a\widetilde{\pi}}=\ol{a}s_\carlitz$.

Comparing with (\ref{exponential}) we also point out, for further references in this paper, 
that (\ref{residues}) becomes in this case:
\begin{equation}\label{residuescarlitz}
s_\carlitz(t)=\sum_{n=0}^\infty\frac{\widetilde{\pi}^{q^n}}{d_n(\theta^{q^n}-t)},\quad |t|<q.
\end{equation}

\subsection{Anderson's functions for elliptic Drinfeld modules\label{edm}}

We recall and deepen some tools described in \cite[Section 4.2.5]{Bourbaki} (see also \cite{ChaPa,Pa}).
Let $z$ be in $\Omega$, and consider the $A$-lattice $\Lambda=\Lambda_z=A+zA$ of rank $2$, 
with associated exponential function $e_z=e_\Lambda$. Let us consider the Drinfeld module $\phi_z$ defined by
\begin{equation}\label{drinfeld_module}
\phi_z:\theta\mapsto \phi_z(\theta)=\theta\tau^0+\tg(z)\tau^1+\tDelta(z)\tau^2,
\end{equation} where
$\tg(z)=\widetilde{\pi}^{q-1}g(z)$, $\tDelta(z)=\widetilde{\pi}^{q^2-1}\Delta(z)$,
with $\Delta=-h^{q-1}$. Then, 
\begin{equation}\label{exponential_functional}\phi_z(a)e_z(\zeta)=e_z(a\zeta)
\end{equation} for all $a\in A$
and $\zeta\in C$ (\cite[Section 5]{Ge}, \cite[Section 4.2.5]{Bourbaki}, see also \cite{Pa}).

We can write, for $\zeta\in C$, 
\begin{equation}\label{ellipticexponential}
e_{z}(\zeta)=\sum_{i=0}^\infty\alpha_i(z)\zeta^{q^i},
\end{equation}
for
functions $\alpha_i:\Omega\rightarrow C$ with $\alpha_0=1$.
By (\ref{exponential_functional}) we deduce, with the initial values $\alpha_0=1,\alpha_{-1}=0$, the recursive 
relations
\begin{equation}\label{alpha_i}
\alpha_i=\frac{1}{[i]}(\widetilde{g}\alpha_{i-1}^q+\widetilde{\Delta}\alpha_{i-2}^{q^2}),\quad i>0.
\end{equation}
This implies that the function $\alpha_i(z)$ is a modular form of weight $q^i-1$ and type $0$ for all $i\geq 0$.
There exist elements $c_{i,m}\in C$ such that
\begin{equation}\label{uexpansionofalpah}\alpha_i(z)=\sum_{m\geq 0}c_{i,m}u^m,\quad i\geq 0,\end{equation} with convergence for $z\in\Omega$ such that $|u|$ is small 
enough. The following lemma tells that a non-zero disk of convergence can be chosen 
independently on $i$.
\begin{Lemme}
We have, for some $B>0$,
\begin{equation}\label{ci0}
c_{i,0}=\frac{1}{d_i}\widetilde{\pi}^{q^i-1},\quad i\geq 0,
\end{equation}
and
\begin{equation}\label{cim}
|c_{i,m}|\leq q^{-q^{i}}B^m,\quad (i,m\geq 0).
\end{equation}
\end{Lemme}

\noindent\emph{Proof.}
Let us write $\widetilde{g}=\sum_{i\geq 0}\widetilde{\gamma}_iu^i$ and $\widetilde{\Delta}=\sum_{i\geq 0}\widetilde{\delta}_iu^i$ with $\widetilde{\gamma}_i,\widetilde{\delta}_i\in C$.
The recursive relations (\ref{alpha_i}) imply, for $i>1$, $m\geq 0$ and $j,k$ non-negative integers:
\begin{equation*}
c_{i,m}=\frac{1}{[i]}\left(\sum_{j+qk=m}\widetilde{\gamma}_jc_{i-1,k}^q+\sum_{j+q^2k=m}\widetilde{\delta}_jc_{i-2,k}^{q^2}\right),
\end{equation*}
from which we deduce at once (\ref{ci0}) because $\widetilde{\gamma}_0=\widetilde{\pi}^{q-1}$ and $\widetilde{\delta}_0=0$.

We now need to provide upper bounds for the $|c_{i,m}|$'s, with explicit dependence on $i,m$. 

Looking at \cite[Definition (5.7), (iii)]{Ge},
there exists $B\geq q$ such that, for all $i\geq 0$, $\max\{|\widetilde{\gamma}_i|,|\widetilde{\delta}_i|\}\leq B^i$. We know that $\alpha_0=1$ and that $|c_{1,m}|\leq q^{-q}B^m$.
After induction and the equality $|[i]|=q^{q^i}$ ($i>0$), we deduce (\ref{cim}).\CVD

In all the following, we shall write:
$$\bfs_1(z,t)=s_{\Lambda_z,z}(t),\quad \bfs_2(z,t)=s_{\Lambda_z,1}(t).$$ These are
functions $\Omega\times B_q\rightarrow C$, where, for $r>0$, $B_r$ is the set $\{t\in C,|t|<r\}$. 

In fact the definition of the functions $s_{\Lambda,\omega}$ tells that
$\bfs_1,\bfs_2\in\mathbf{Hol}(\Omega)[[t]]$, where $\mathbf{Hol}(\Omega)$ denotes the $C$-algebra of rigid holomorphic functions $\Omega\rightarrow C$.
After (\ref{ellipticexponential}) and (\ref{residues}) we see that, for any couple $(z,t)\in\Omega\times B_q$, the following convergent series expansions hold:
\begin{eqnarray*}
\bfs_1(z,t)&=&\sum_{i=0}^\infty\frac{\alpha_i(z)z^{q^i}}{\theta^{q^i}-t}\\
\bfs_2(z,t)&=&\sum_{i=0}^\infty\frac{\alpha_i(z)}{\theta^{q^i}-t}.
\end{eqnarray*}

Our notations stress the dependence on two variables $z\in\Omega,t\in B_q$. For these functions, we will also write, occasionally,
$\bfs_1(z),\bfs_2(z)$, to stress the dependence on $z\in\Omega$.
We can also fix $z\in\Omega$ and study the functions $\bfs_1(z,\cdot),\bfs_2(z,\cdot)$, or look at the functions $\bfs_1(\cdot,t),\bfs_2(\cdot,t):\Omega\rightarrow\TT_{<q}$
with formal series as values. In the next section, we provide the necessary analysis of the functions $\bfs_1(z,\cdot),\bfs_2(z,\cdot)$. Hence, we fix now $z\in\Omega$.

\subsubsection{The $\bsb{s}_i$'s as functions of the variable $t$, with $z$ fixed.\label{functionsofthevariable}}

At $\theta$, the functions $\bfs_i(z,\cdot)$ have simple poles. Their respective residues are, according to
Section \ref{section:trivialisations}, $-z$ for the function $\bfs_1(z,\cdot)$ and $-1$
for $\bfs_2(z,\cdot)$.
Moreover, we have $\bfs_1^{(1)}(z,\theta)=\eta_1$ and $\bfs_2^{(1)}(z,\theta)=\eta_2$, where $\eta_1,\eta_2$
are the {\em quasi-periods} of $\Lambda_z$ (see \cite[Section 4.2.4]{Bourbaki} and \cite[Section 7]{gekeler:compositio}).

Let us consider the matrix function:
$$\hPsi (z,t):=\sqm{\bfs_1(z,t)}{\bfs_2(z,t)}{\bfs^{(1)}_1(z,t)}{\bfs_2^{(1)}(z,t)}.$$
By \cite[Section 4.2.3]{Bourbaki} (see in particular equation (15)), we have:
\begin{equation}\label{fundamentalrelation}
\hPsi(z,t)^{(1)}=\widetilde{\Theta}(z)\cdot\hPsi(z,t),\quad \text{where }\widetilde{\Theta}(z)=\sqm{0}{1}{\frac{t-\theta}{\tDelta(z)}}{-\frac{\tg(z)}{\tDelta(z)}},
\end{equation} 
yielding the following $\tau$-difference linear equation of order $2$:
\begin{equation}\label{fundamentalrelation2}
\bfs_2^{(2)}=\frac{t-\theta}{\tDelta}\bfs_2-\frac{\tg}{\tDelta}\bfs_2^{(1)}.
\end{equation}

\begin{Remarque}{\em 
By \cite{An}, there is a fully faithful contravariant functor from the 
category of Drinfeld $A$-modules over $K^{\text{alg.}}$ to the category of {\em Anderson's $A$-motives} over $K^{\text{alg.}}$. 
Part of this association is sketched in \cite[Section 4.2.2]{Bourbaki}, where the definition of $A$-motive is given and discussed
(see also \cite{ChaPa}); it is based precisely
on Anderson's functions $s_{\Lambda,\omega}$. In the language introduced by Anderson, $\hPsi$ is a {\em rigid analytic trivialisation} of the 
$A$-motive associated to the Drinfeld module $\phi=\phi_\Lambda$.}\end{Remarque}

We will also use the following fundamental lemma, whose proof depends on Gekeler's paper \cite{gekeler:compositio}.

\begin{Lemme}[``Deformation of Legendre's identity"]\label{lemmalegendre}
We have, for all $z\in\Omega$ and $t$ with $|t|<q$:
\begin{equation}\label{dethPsi}
\det(\hPsi)=\widetilde{\pi}^{-1-q}h(z)^{-1}s_\carlitz(t).
\end{equation}
\end{Lemme}

\noindent\emph{Proof.} Let $f(z,t)$ be the function $\det(\hPsi(z,t))h(z)\widetilde{\pi}^{1+q}$, for $z\in\Omega$ and $t\in B_q$. 
We have:
$$f^{(1)}(z,t)=-(t-\theta)\tDelta(z)^{-1}\det(\hPsi(z,t))h(z)^q\widetilde{\pi}^{q+q^2}=(t-\theta)f(z,t).$$ 

For fixed $z\in\Omega$, we know that
$\bfs_i^{(k)}(z,\cdot)\in\TT_{<q^{q^{k}}}\subset\TT$ for all $k\geq 0$. Hence, $f(z,\cdot)\in\TT$ for all $z\in\Omega$. By arguments used in
the remark on the $\overline{K}$-vector space structure of the set of solutions 
of (\ref{taudifferencebfs}), $f(z,t)$ is equal to $\lambda(z,t)s_\carlitz(t)$, for some $\lambda(z,t)\in \overline{A}$; the matter is now to compute $\lambda$, which does not depend on $z\in\Omega$ as  it follows easily by fixing $t=t_0\in B_q$ transcendental over $\FF_q$ 
and observing that $f(z,t_0)$ is holomorphic over $\Omega$ with values in a discrete
set.

Now, for $z$ fixed as $t\rightarrow\theta$, 
$$\lim_{t\rightarrow\theta}\hPsi(z,t)-\begin{pmatrix}-\frac{z}{t-\theta} & -\frac{1}{t-\theta}\\ \eta_1& \eta_2\end{pmatrix}=\begin{pmatrix} * & * \\ 0 & 0\end{pmatrix},$$
$\eta_1,\eta_2$ being the 
{\em quasi-periods} (periods of second kind) of the lattice $A\omega_1+A\omega_2$ (respectively associated to $\omega_1$ and $\omega_2$)
\cite[Section 7, Equations (7.1)]{gekeler:compositio}, with generators $\omega_1=z,\omega_2=1$, where the asterisks denote continuous functions of the variable $z$.
Hence, we have $\lim_{t\rightarrow\theta}(t-\theta)\det(\hPsi(z,t))=-z\eta_2+\eta_1$.
By \cite[Theorem 6.2]{gekeler:compositio}, $-z\eta_2+\eta_1=-\widetilde{\pi}^{-q}h(z)^{-1}$.

At once:
\begin{eqnarray*}
\lefteqn{-\widetilde{\pi}^{-q}h(z)^{-1}=}\\
&=&\lim_{t\rightarrow\theta}(t-\theta)\det(\hPsi(z,t))\\
&=&\lambda(\theta)\widetilde{\pi}^{-q-1}h(z)^{-1}\lim_{t\rightarrow\theta}(t-\theta)s_\carlitz(t)\\
&=&-\lambda(\theta)\widetilde{\pi}^{-q}h(z)^{-1},
\end{eqnarray*}
which implies that $\lambda=\lambda(\theta)=1$ ($\theta$ is transcendental over $\FF_q$). Our Lemma follows.\CVD

In the next section, we study the functions $\bsb{s}_1,\bsb{s}_2$ as functions $\Omega\rightarrow\TT_{<q}$.

\subsubsection{The $\bsb{s}_i$'s as functions $\Omega\rightarrow\TT_{<q}$.\label{modularity}}

We observe, by the definitions of $\bfs_1,\bfs_2$, and by the fact, remarked in (\ref{alpha_i}), that $\alpha_i$ is a modular form of weight $q^i-1$ and type $0$ for all $i$,
and by (\ref{pseudolinearity}), that for all $\gamma=\sqm{a}{b}{c}{d}\in\Gamma$:
\begin{eqnarray*}
\bfs_2(\gamma(z),t)&=&\sum_{i=0}^\infty(cz+d)^{q^i-1}\frac{\alpha_i(z)}{\theta^{q^i}-t}\\
&=&(cz+d)^{-1}s_{\Lambda_z,cz+d}(z)\\
&=&(cz+d)^{-1}(\overline{c}\bfs_1(z,t)+\overline{d}\bfs_2(z,t)).
\end{eqnarray*}
Similarly,
\begin{eqnarray*}
\bfs_1(\gamma(z),t)&=&\sum_{i=0}^\infty(cz+d)^{q^i-1}\frac{\alpha_i(z)(\gamma(z))^{q^i}}{\theta^{q^i}-t}\\
&=&(cz+d)^{-1}s_{\Lambda_z,az+b}(z)\\
&=&(cz+d)^{-1}(\overline{a}\bfs_1(z,t)+\overline{b}\bfs_2(z,t)).
\end{eqnarray*}
Let us write 
$$\Sigma(z,t):=\binom{\bfs_1(z,t)}{\bfs_2(z,t)}.$$ We have proved:
\begin{Lemme}\label{lemmavectorial} For all $\gamma=\sqm{a}{b}{c}{d}\in\Gamma$, and for all $z\in\Omega$, we have the following identity of series in $\TT_{<q}$:
\begin{equation}\label{eqvectorial}
\Sigma(\gamma(z),t)=(cz+d)^{-1}\overline{\gamma}\cdot\Sigma(z,t),\end{equation} 
where $\overline{\gamma}$ is the matrix $\sqm{\overline{a}}{\overline{b}}{\overline{c}}{\overline{d}}\in\Gamma$.
\end{Lemme}

\subsubsection{Behaviour of $\bsb{s}_2$ at the infinity cusp and $u$-expansion\label{behaviouru}}

We use the results of the previous subsections to see how the function $\bfs_2$ behaves for $z$ approaching the cusp at infinity
of the rigid analytic space $\Gamma\backslash\Omega$. 
Here we will prove two  lemmas.

\begin{Lemme}\label{usefulford}
There exists a real number $r>0$ such that for all $(z,t)\in\Omega\times C$ with $|u|=|u(z)|<r,|t|<r$, we have:
\begin{equation}
\bfs_2(z,t)=\widetilde{\pi}^{-1}s_\carlitz(t)+\sum_{m\geq 1}\kappa_m(t)u^m,\label{uexpansions2}
\end{equation}
where
for $m\geq 1$,
\begin{eqnarray*}
\kappa_m(t)&=&\sum_{i\geq 1}\frac{c_{i,m}}{\theta^{q^i}-t}\\
&=&\sum_{j\geq 0}t^j\sum_{i\geq 1}c_{i,m}\theta^{-q^i(1+j)}\in \TT_{<q^q},
\end{eqnarray*}
the $c_{i,m}$'s being the coefficients in the expansions (\ref{uexpansionofalpah}).
\end{Lemme}

\noindent\emph{Proof.} For $z\in\Omega$ such that $|u|<B^{-1}$ with $B$ as in (\ref{cim}), and for $|t|<q$, (\ref{residues}) yields:
\begin{eqnarray*}
\bfs_2(z,t)&=&\frac{1}{\theta-t}+\sum_{i\geq 1}\frac{\alpha_i(z)}{\theta^{q^i}-t}\nonumber\\
&=&\frac{1}{\theta-t}+\sum_{i\geq1}\sum_{m\geq0}c_{i,m}u^m\frac{1}{\theta^{q^i}-t}\nonumber\\
&=&\frac{1}{\theta-t}+\sum_{i\geq 1}c_{i,0}\frac{1}{\theta^{q^i}-t}+\sum_{m\geq 1}u^m\sum_{i\geq 1}\frac{c_{i,m}}{\theta^{q^i}-t}\nonumber\\
&=&\widetilde{\pi}^{-1}\sum_{i\geq 0}\frac{\widetilde{\pi}^{q^i}}{d_i}\frac{1}{\theta^{q^i}-t}+\sum_{m\geq 1}\kappa_m(t)u^m\nonumber\\
&=&\widetilde{\pi}^{-1}s_\carlitz(t)+\sum_{m\geq 1}\kappa_m(t)u^m,
\end{eqnarray*}
where, in the second equality we have substituted the $u$-expansions of the $\alpha_i$'s in our formulas, in the third we have separately considered constant terms, 
in the fourth equality, we have used (\ref{ci0}), in the fifth we have recognised the shape of $s_\carlitz$ 
(\ref{residuescarlitz}),
and we have noticed, by using (\ref{cim}), that for all $t\in C$ such that
$|t|\leq q$, $|\kappa_m(t)|\leq B^mq^{-1}$. \CVD

Later, we will need to do some arithmetic with the $u$-expansion (\ref{uexpansions2}).
To this purpose, it is advantageous to set:
$$\bsb{d}(z,t):=\widetilde{\pi}s_\carlitz(t)^{-1}\bfs_2(z,t),$$ function for which (\ref{fundamentalrelation2}) becomes:
\begin{equation}\label{equexpansion}
\bsb{d}=(t-\theta^q)\Delta\bsb{d}^{(2)}+g\bsb{d}^{(1)}.
\end{equation}
We will need part of the following lemma.

\begin{Lemme}\label{lemmed}
We have 
\begin{equation}
\bsb{d}=\sum_{i\geq 0}c_i(t)u^{(q-1)i}\in1+u^{q-1}\FF_q[t,\theta][[u^{q-1}]].
\end{equation}
More precisely, 
$$\bsb{d}=1+(\theta-t)u^{q(q-1)}+(\theta-t)u^{(q^2-q+1)(q-1)}+\cdots\in1+(t-\theta)u^{q-1}\FF_q[t,\theta][[u^{q-1}]],$$
where the dots $\cdots$ stand for terms of higher order in $u$.

Let $i$ be a positive integer. We have
$$-\infty\leq\deg_tc_i\leq \log_{q^2}i,$$ where $\log_{q^2}$ is the logarithm in base $q^2$, with the convention $\deg_t0=-\infty$.

\end{Lemme}

\noindent\emph{Proof.} For simplicity, we write $v=u^{q-1}$. It is clear, looking at Lemma \ref{usefulford}, that $\bsb{d}$ is a series in $\TT_{<q^q}[[v]]$.
We have the series expansions (cf. \cite[Section 10]{Ge}):
\begin{eqnarray*}
g&=&1-[1]v+\cdots=\sum_{n=0}^\infty \gamma_{n}v^n\in A[[v]],\\
\Delta&=&-v(1-v^{q-1}+\cdots)=\sum_{n=0}^\infty \delta_{n}v^n\in uA[[v]],
\end{eqnarray*}
We deduce, from (\ref{equexpansion}), that 
\begin{equation}\label{equationci}
c_m=(t-\theta^q)\sum_{i+q^2j=m}\delta_ic_j^{(2)}+\sum_{i+qj=m}\gamma_ic_j^{(1)},
\end{equation}
which yields inductively that $c_i$ belongs to $\FF_q[t,\theta]$, because the coefficients of the $u$-expansions
of $\Delta$ and $g$ are $A$-integral. The statement on the degrees of the coefficients
of the $c_i$'s, is also a simple inductive consequence of (\ref{equationci}) and the following two facts: that $\deg_t\delta_i,\deg_t\gamma_i\leq 0$,
and that $\deg_tc_i^{(k)}=\deg_tc_i$ for all $i,k$ ($t$ is $\tau$-invariant).

The explicit formula for the coefficients $c_i$ with $i\leq q^2-q+1$ is an exercice that we leave to the reader, which needs 
\cite[Corollaries (10.3), (10.11)]{Ge}. The explicit computation can be pushed easily to coefficients of higher order, but we 
skip it as we will not need these explicit formulas at all in this paper. The fact that the coefficients $c_i$ belong to the ideal generated by $t-\theta$
for $i\geq 1$ follows from the computation of the residues in \ref{functionsofthevariable}.
\CVD

\section{The function $\bsb{E}$\label{bsbE}}

The function of the title is defined, for $z\in\Omega$ and $t\in B_q$, by:
$$\boldsymbol{E}(z,t)=-h(z)\bsb{d}^{(1)}(z,t)=-(t-\theta)^{-1}\widetilde{\pi}^qh(z)s_\carlitz^{-1}(t)\boldsymbol{s}_2^{(1)}(z,t),$$ with $\bsb{d}$ the function
of Lemma \ref{lemmed}.
This section is entirely devoted to the description of its main properties. Three Propositions will be proved here.

In Proposition \ref{taudifferenceeq} we use the arguments of \ref{functionsofthevariable} to show that, just as $\bsb{d}$, $\bsb{E}$ satisfies a linear $\tau$-difference equation of order $2$ with coefficients isobaric in $C[[t]][g,h]$ (\footnote{This phenomenon holds with more generality and should be compared with
a result of Stiller in \cite{Stiller}.}).

In Proposition \ref{propositionautomorphyE}, where we use this time the arguments developed in \ref{modularity}, we analyse the functional equations relating the values of $\bsb{E}$ at $z$ and $\gamma(z)$,
where $\gamma=\begin{pmatrix}a&b\\ c&d\end{pmatrix}\in\Gamma$; they involve the {\em factors of automorphy}:
$$J_\gamma(z)=cz+d,\quad \bsb{J}_\gamma(z)=\ol{c}\frac{\bsb{s}_1(z,t)}{\bsb{s}_2(z,t)}+\ol{d},$$ with values convergent in $C[[t]]$.

Proposition \ref{propositionEuexpansion} follows from what we did in \ref{behaviouru} and describes the third important feature of the function $\bsb{E}$; the existence of a $u$-expansion
in $\FF_q[t,\theta][[u]]$. For Drinfeld quasi-modular forms, the degree in $\theta$ of the $n$-th coefficient of the $u$-expansion grows pretty rapidly with $n$
in contrast of the classical framework. The function $\bsb{E}$ does not make exception to this principle. However, the degree in $t$ of the $n$-th coefficient 
grows slowly, and this property is used crucially in the proof of the multiplicity estimate. Another important property studied in this section is that $\bsb{E}(z,\theta)$ is
a well defined function $\Omega\rightarrow C$ and is equal to Gekeler's
function $E$.

\subsection{linear $\tau$-difference equations}

\begin{Proposition}\label{taudifferenceeq} For all $z\in\Omega$, the function $\bsb{E}(z,\cdot)$
can be developed as a series of $\TT_{< q^q}$. Moreover,
The following linear $\tau$-difference equation holds  in $\TT_{<q^q}$, for all $z\in\Omega$:
\begin{equation}\label{eq:sigmaboldE}
\boldsymbol{E}^{(2)}=\frac{1}{t-\theta^{q^2}}(\Delta\boldsymbol{E}+g^q\boldsymbol{E}^{(1)}).
\end{equation}
\end{Proposition}
\noindent\emph{Proof.} After having chosen a $(q-1)$-th root of $-\theta$,
let us write, following Anderson, Brownawell, and Papanikolas in \cite[Section 3.1.2]{ABP}, $$\boldsymbol{\Omega}(t):=(-\theta)^{\frac{-q}{q-1}}\prod_{n=1}^\infty\left(1-\frac{t}{\theta^{q^n}}\right)\in (\TT_\infty\cap K_\infty((-\theta)^{\frac{1}{q-1}})[[t]])\setminus K_\infty(t)^{\text{alg}}.$$
It is plain that
$$\boldsymbol{\Omega}^{(-1)}(t)=(t-\theta)\boldsymbol{\Omega}(t).$$
Thanks to the remark on the $\overline{K}$-vector space structure of the set of solutions 
of (\ref{taudifferencebfs}) and after the computation of the constant of proportionality, we get
\begin{equation}\label{omegas}
s_\carlitz(t)=\frac{1}{\boldsymbol{\Omega}^{(-1)}(t)}.
\end{equation}
At once, we obtain that the function $s_\carlitz$ has no zeros in the domain $C\setminus\{\theta,\theta^q,\ldots\}$ from which it follows that $((t-\theta)s_\carlitz)^{-1}\in\TT_{<q^q}$. Moreover, for all $z\in\Omega$, we have $\bfs_2\in\TT_{<q}$ so that
$\bfs_2^{(1)}\in\TT_{<q^q}$. Multiplying the factors that define the function $\boldsymbol{E}$, we then get, for all
$z\in\Omega$, that $\boldsymbol{E}(z,\cdot)\in\TT_{<q^q}$, which gives the first part of the proposition (and in fact, it can be proved that $\bsb{d},\bsb{E}(z,\cdot)\in\TT_\infty$
for all $z\in\Omega$, but we skip on this property since it will not be needed in the present paper).

In order to prove the second part of the proposition, we remark, 
from (\ref{equexpansion})  (or what is the same, (\ref{fundamentalrelation2})), that
$$\bfs_2^{(3)}=\frac{t-\theta^q}{\widetilde{\Delta}^q}\bfs_2^{(1)}-\frac{\widetilde{g}^q}{\widetilde{\Delta}^q}\bfs_2^{(2)},\quad\text{ or equivalently, }\quad
\bsb{d}^{(3)}=\frac{1}{(t-\theta^{q^2})\Delta^q}(\bsb{d}^{(1)}-g^q\bsb{d}^{(2)}).$$
By the definition of $\boldsymbol{E}$ and the $\tau$-difference equation (\ref{taudifferencebfs}) we find the relation:
\begin{eqnarray}
\boldsymbol{E}^{(k)}&=&-(t-\theta^{q^k})^{-1}(t-\theta^{q^{k-1}})^{-1}\cdots(t-\theta)^{-1}\widetilde{\pi}^{q^{k+1}}h^{q^k}s_{\carlitz}^{-1}\bfs_2^{(k+1)},\nonumber\\
&=&-h^{q^k}\bsb{d}^{(k+1)}\label{Eintermsofbfs2}\end{eqnarray} for $k\geq 0$.
Substituting the above expression for $\bsb{d}^{(3)}$ in it, we get what we expected.\CVD

\subsection{Factors of automorphy, modularity\label{JL}}

In the next proposition, the function $\bsb{E}$ is viewed as a function $\Omega\rightarrow\TT_{>0}$ (it can be proved that it defines, in fact, a function
$\Omega\rightarrow\TT_\infty$).
In order to state the proposition, we first need a preliminary discussion.

If $\omega\not\in\theta\Lambda$, then
$e_\Lambda(\omega/\theta)\neq 0$ and $s_{\Lambda,\omega}(t)\in\TT_{>0}^\times$ (group of units of $\TT_{>0}$), so that, for every $z$ fixed, $\bfs_2(z,\cdot)^{-1}\in\TT_{>0}^\times$ (\footnote{The radius of convergence, in principle depending on $z$, seems difficult to compute.}). Hence, we have a well defined map
$$\begin{array}{rclcl}\bsb{\xi} &:&\Omega & \rightarrow & \TT_{>0}^\times\\
& & z & \mapsto& \frac{\bfs_1(z,t)}{\bfs_2(z,t)},\end{array}$$
and we can consider the map
$$(\gamma,z)=\left(\begin{pmatrix}a&b\\ c&d\end{pmatrix},z\right)\in\Gamma\times\Omega\mapsto \bsb{J}_\gamma(z):=\ol{c}\bsb{\xi}+\ol{d}\in\TT_{>0}.$$
Since $c,d$ are relatively prime, we have $cz+d\not\in\theta\Lambda_z$ implying that $\ol{c}\bsb{s}_1+\ol{d}\bsb{s}_2=s_{\Lambda_z,cz+d}\in\TT_{>0}^\times$.
Therefore, for all $\gamma\in\Gamma$ and $z\in\Omega$, $\bsb{J}_\gamma\in\TT_{>0}^\times$.

Moreover, by (\ref{eqvectorial}) we have, for all $\gamma\in\Gamma$ and $z\in\Omega$, 
\begin{equation}\label{eq:iotazeta}
\bsb{\xi}(\gamma(z))=\overline{\gamma}(\bsb{\xi}(z))\in C((t)),
\end{equation}
so that, for $\gamma,\delta\in\Gamma$ and $z\in\Omega$, 
\begin{equation}\bsb{J}_{\gamma\delta}(z)=\bsb{J}_\gamma(\delta(z))\bsb{J}_{\delta}(z).
\label{eq:automorphy2}\end{equation}

the map $\bsb{J}:\Gamma\times\Omega\rightarrow\TT_{>0}^\times$ is our ``new" {\em factor of automorphy}, to be considered together with the more familiar factor of automorphy
$$J_\gamma(z):=cz+d.$$

Let us also write, for $\gamma=\sqm{a}{b}{c}{d}\in\Gamma$:

\begin{eqnarray*}
L_\gamma(z)&=&\frac{c}{cz+d},\\
\boldsymbol{L}_\gamma(z)&=&\frac{\overline{c}}{\overline{c}\bfs_1+\overline{d}\bfs_2}.\end{eqnarray*}
We remark that for all $\gamma\in\Gamma$, $\boldsymbol{L}_\gamma(z)$ belongs to $\TT_{>0}$ because $\bsb{s}_2\bsb{J}_\gamma(z)\in\TT_{>0}^\times$.
Moreover, the functions $\bsb{J}_\gamma$ and $(\theta-t)^{-1}\bsb{L}_\gamma$ are deformations of
$J_\gamma$ and $L_\gamma$ respectively, for all $\gamma\in\Gamma$.
Indeed, we recall that $(t-\theta)\bfs_2(z,t)\rightarrow-1$ and $(t-\theta)\bfs_1(z,t)\rightarrow -z$ as $t\rightarrow \theta$. 
Hence, $\lim_{t\rightarrow\theta}\frac{\bsb{s}_1}{\bsb{s}_2}=z$.
This implies that \begin{equation}\label{Jbold}\lim_{t\rightarrow\theta}\boldsymbol{J}_\gamma(z)=J_\gamma(z).\end{equation}
In a similar way we see that 
\begin{equation}\label{Lbold}\lim_{t\rightarrow\theta}(t-\theta)^{-1}\boldsymbol{L}_\gamma(z)= -L_\gamma(z).\end{equation}

We further define the sequence of 
functions
$(g^\star_k)_{k\geq 0}$ by:
\[g_{-1}^\star=0,\quad g_0^\star=1,\quad g_1^\star=g,\quad g_k^\star=(t-\theta^{q^{k-1}})g_{k-2}^\star \Delta^{q^{k-2}}+g_{k-1}^\star g^{q^{k-1}},\quad k\geq 2,\]
so that for all $k\geq 0$, we have the identity $g^\star_k(z,\theta)=g_k(z)$, the function introduced in \cite[Equation (6.8)]{Ge}.

We have:
\begin{Proposition}\label{propositionautomorphyE}
For all $z\in\Omega$, $\gamma\in\Gamma$ and $k\geq 0$ the following identity of formal series of $\TT_{>0}$ holds:
\begin{eqnarray}
\boldsymbol{E}^{(k)}(\gamma(z),t)&=&\det(\gamma)^{-1}J_\gamma(z)^{q^k}\bsb{J}_\gamma(z)\times\label{tmodularityE}\\
& &\left(\boldsymbol{E}^{(k)}(z,t)+\frac{g_k^\star(z)}{\widetilde{\pi} (t-\theta)(t-\theta^q)\cdots(t-\theta^{q^k})}\bsb{L}_\gamma(z)\right).\nonumber\end{eqnarray}
\end{Proposition}

\noindent\emph{Proof.} From the deformation of Legendre's identity (\ref{dethPsi}) we deduce that
\begin{equation}\label{bfs1gamma}
\bfs_1^{(1)}=\frac{1}{\bfs_2}(\bfs_1\bfs_2^{(1)}-\widetilde{\pi}^{-1-q}h^{-1}s_\carlitz).
\end{equation}
Let $\gamma=\begin{pmatrix}a & b \\ c & d\end{pmatrix}\in\Gamma$.
Applying $\tau$ on both left and right hand sides of
\begin{equation}\label{bfs0tra}
\bfs_2(\gamma(z))=J_\gamma^{-1}\boldsymbol{J}_\gamma\bfs_2(z)=J_\gamma^{-1}(\ol{c}\bfs_1(z)+\ol{d}\bfs_2(z)),
\end{equation}
consequence of Lemma \ref{lemmavectorial}, we see that
$$\bfs_2^{(1)}(\gamma(z))=J_\gamma^{-q}(\ol{c}\bfs_1^{(1)}+\ol{d}\bfs_2^{(1)}).$$
We now eliminate $\bsb{s}_1^{(1)}$ from this identity and (\ref{bfs1gamma}), getting identities in $\TT_{>0}$. Indeed,
\begin{eqnarray*}
\lefteqn{\ol{c}\bfs_1^{(1)}+\ol{d}\bfs_2^{(1)}=}\\
&=&\frac{\ol{c}}{\bfs_2}(\bfs_1\bfs_2^{(1)}-\widetilde{\pi}^{-1-q}h(z)^{-1}s_\carlitz(t))+\ol{d}\bfs_2^{(1)}\\
&=&\bfs_2^{(1)}\left(\ol{c}\frac{\bfs_1}{\bfs_2}+\ol{d}\right)-\widetilde{\pi}^{-1-q}h(z)^{-1}s_\carlitz(t)\bfs_2^{-1}\\
&=&\left(\ol{c}\frac{\bfs_1}{\bfs_2}+\ol{d}\right)\left(\bfs_2^{(1)}-\widetilde{\pi}^{-1-q}h(z)^{-1}s_\carlitz(t)\frac{\ol{c}}{\ol{c}\bfs_1+\ol{d}\bfs_2}\right),
\end{eqnarray*}
that is,
\begin{equation}\label{functionalequationstwist}
\bfs_2^{(1)}(\gamma(z))=J_\gamma^{-q}\boldsymbol{J}_\gamma\left(\bfs_2^{(1)}(z)-\frac{\widetilde{\pi}^{-1-q}s_\carlitz(t)}{h(z)}\boldsymbol{L}_\gamma\right).
\end{equation}
This functional equation is equivalent to the following functional equation for $\bsb{d}^{(1)}$ (in $\TT_{>0}$):
\begin{equation}\label{bfs1tra}
\bsb{d}^{(1)}(\gamma(z))=J_\gamma^{-q}\boldsymbol{J}_\gamma\left(\bsb{d}^{(1)}(z)-\frac{1}{\widetilde{\pi}(t-\theta)h(z)}\boldsymbol{L}_\gamma\right).
\end{equation}
This already implies, by the definition of $\bsb{E}$ and the modularity of $h$:
\begin{equation*}
\boldsymbol{E}(\gamma(z))=\det(\gamma)^{-1}J_\gamma\boldsymbol{J}_\gamma\left(\boldsymbol{E}(z)+\frac{1}{\widetilde{\pi}(t-\theta)}\boldsymbol{L}_\gamma\right)
\end{equation*}
which is our proposition for $k=0$.

We point out that (\ref{eqvectorial}) implies the functional equation, for all $\gamma\in\Gamma$:
\begin{equation}\label{kequaltominusone}
\bsb{d}(\gamma(z))=J_\gamma^{-1}\bsb{J}_\gamma\bsb{d}(z).
\end{equation}
The joint application of (\ref{kequaltominusone}), (\ref{bfs0tra}), (\ref{bfs1tra}) and (\ref{fundamentalrelation}) and induction on $k$ imply, for all $k\geq 0$ and 
$\gamma\in\Gamma$, the functional equation in $\TT_{>0}$:
\begin{equation}
\bsb{d}^{(k)}(\gamma(z))=J_\gamma^{-q^k}\boldsymbol{J}_\gamma\left(\bsb{d}^{(k)}(z)-\frac{g^\star_{k-1}}{\widetilde{\pi}h(z)^{q^{k-1}}(t-\theta)(t-\theta^q)\cdots(t-\theta^{q^{k-1}})}\bsb{L}_\gamma\right),
\end{equation}
where we have also used the functional equation (\ref{taudifferencebfs}).
By (\ref{Eintermsofbfs2}),
we end the proof of the proposition.\CVD

\subsection{$u$-expansions}


\begin{Proposition}\label{propositionEuexpansion}
We have $$\boldsymbol{E}(z,t)=u\sum_{n\geq 0}c_n(t)u^{(q-1)n}\in u\FF_q[\theta,t][[\unif^{q-1}]],$$
where the formal series on the right-hand side converges 
for all $t,u$ with $|t|\leq q$ and $|u|$ small. The terms of order $\leq q(q-1)$ of the $u$-expansion of $\boldsymbol{E}$
are:
\begin{equation}\label{eq:u-expansionofE}
\boldsymbol{E}=u(1+u^{(q-1)^2}-(t-\theta)u^{(q-1)q}+\cdots).\end{equation}
Moreover, for all $n>0$, we have the following inequality for the degree in $t$
of $c_n(t)$:
$$\deg_tc_n\leq\log_q n,$$ where $\log_q$ denotes the logarithm in base $q$ and 
where we have adopted the convention $\deg_t0=-\infty$.
\end{Proposition}

\noindent\emph{Proof.} This is a simple consequence of Lemma \ref{lemmed} and the definition of $\bsb{E}$.\CVD


\begin{Remarque}{\em 
Let us introduce the function $$\boldsymbol{\mu}=\widetilde{\pi}^{1-q}\bfs_2^{(1)}/\bfs_2\in C[[t,\unif^{q-1}]].$$
By (\ref{fundamentalrelation}), $\bsb{\mu}$ satisfies the non-linear $\tau$-difference equation:
$$\boldsymbol{\mu}^{(1)}=\frac{(t-\theta)}{\Delta}\boldsymbol{\mu}^{-1}-\frac{ g }{\Delta}.$$
Hence, $\boldsymbol{\mu}=(t-\theta)\Delta^{-1}(\boldsymbol{\mu}^{(1)}+ g /\Delta)^{-1}$. Although not needed 
in this paper, we point out that this functional equation gives the following 
continued fraction development, which turns out to be convergent for the $\unif$-adic topology:
\begin{equation}\label{continuedfractionmu}\boldsymbol{\mu}=\cfrac{(t-\theta)}{ g +\cfrac{\Delta(t-\theta^q)}{ g^q+\cfrac{\Delta^q(t-\theta^{q^2})}{ g^{q^2}+
\cfrac{\Delta^{q^2}(t-\theta^{q^3})}{\dotsb}}}}\in \FF_q[t,\theta][[v]].\end{equation}
This property should be compared with certain continued fraction developments in \cite[Section 4, 5]{KZ1},
or the continued fraction developments described after \cite[Theorem 2]{KK1}.}
\end{Remarque}

\section{Bi-weighted automorphic functions\label{biweighted}}
In this section we introduce a class of bi-weighted automorphic functions that we call {\em almost $A$-quasi-modular forms}. We will see that they generate a $\TT_{>0}$-algebra $\widetilde{\mathcal{M}}$ with natural embedding in $C[[t,u]]$.
Thanks to
the two kinds of factor of automorphy described below, $\widetilde{\mathcal{M}}$ is also graded by the group $G=\ZZ^2\times\ZZ/(q-1)\ZZ$.
We will not pursue, in this paper, any investigation on the structure of $\widetilde{\mathcal{M}}$; this will be objective of another
work.

We will show, with the help of the results of Section \ref{bsbE}, that 
$g,h,\bsb{E},\bsb{F}\in\widetilde{\mathcal{M}}$ with $\bsb{F}=\tau\bsb{E}$.
It will be proved that for this graduation, the degrees (in $\ZZ^2\times\ZZ/(q-1)\ZZ$) of these functions are respectively the following elements of $G$: $(q-1,0,0),(q+1,1,0),(1,1,1)$ and $(q,1,1)$ and we will show from this that 
they are algebraically independent over $C((t))$ (also $E$ belongs to $\widetilde{\mathcal{M}}$, but we will not use this property).
Since they take values in $\TT_{<q^q}$, we will study with some detail the four dimensional $\TT_{<q^q}$-algebra
$$\MM^\dag:=\TT_{<q^q}[g,h,\bsb{E},\bsb{F}].$$ 

Proposition \ref{taudifferenceeq} implies that $\tau$ acts on $\MM^\dag$:
If $\bsb{f}\in\MM^\dag$ is homogeneous of degree $(\mu,\nu,m)$ then $\tau\bsb{f}$ is also homogeneous of degree $(q\mu,\nu,m)$. 

We will see that  if $\bsb{f}\in\MM^\dag$ is homogeneous of degree $(\mu,\nu,m)$, the function
$$\begin{array}{rrcl}&\Omega&\rightarrow& C \\ \varepsilon(\bsb{f}):&z&\mapsto&\bsb{f}(z)|_{t=\theta}\end{array}$$ is a well defined Drinfeld quasi-modular form of weight $\mu+\nu$, type $m$
and depth $\leq\nu$. An example is given by Lemma \ref{EEgivesE}: $\varepsilon(\bsb{E})=E$.

\subsection{Preliminaries on the functions $\bsb{J}_\gamma$ and $\bsb{L}_\gamma$}
Let us consider three matrices in $\Gamma$:
\begin{equation}\label{ABC}
\mathcal{A}=\begin{pmatrix}a & b \\ c & d\end{pmatrix},\quad \mathcal{B}=\begin{pmatrix}\alpha & \beta \\ \gamma & \delta\end{pmatrix},\quad \mathcal{C}={\mathcal{A}}\cdot{\mathcal{B}}
=\begin{pmatrix}* & * \\ x & y\end{pmatrix}\in\Gamma.
\end{equation}

\begin{Lemme}\label{lemmeABC}
We have the following identities in $\TT_{>0}$:
\begin{eqnarray*}
L_{\mathcal{A}}(\mathcal{B}(z))&=&\det(\mathcal{B})^{-1}J_{\mathcal{B}}(z)^2(L_{\mathcal{C}}(z)-L_{\mathcal{B}}(z)),\\
\boldsymbol{L}_{\mathcal{A}}(\mathcal{B}(z))&=&\det(\mathcal{B})^{-1}J_{\mathcal{B}}(z)\boldsymbol{J}_{\mathcal{B}}(z)(\boldsymbol{L}_{\mathcal{C}}(z)-\boldsymbol{L}_{\mathcal{B}}(z)).
\end{eqnarray*}
\end{Lemme}
\noindent\emph{Proof.} We begin by proving the first formula, observing that $c=\det(\mathcal{B})^{-1}(x\delta-y\gamma)$:
\begin{eqnarray*}
\lefteqn{J_{\mathcal{B}}(z)^2(L_{\mathcal{C}}(z)-L_{\mathcal{B}}(z))=}\\
&=&(\gamma z+\delta)^2\left(\frac{x}{xz+y}-\frac{\gamma}{\gamma z+\delta}\right)\\
&=&\det(\mathcal{B})\frac{c(\delta+\gamma z)}{(c\alpha+d\gamma)z+(c\beta+d\gamma)}\\
&=&\det(\mathcal{B})\frac{c}{\frac{(c\alpha+d\gamma)z+(c\beta+d\gamma)}{\delta+\gamma z}}\\
&=&\det(\mathcal{B})\frac{c}{c\frac{\alpha z+\beta}{\gamma z+\delta}+d}\\
&=&\det(\mathcal{B})L_{\mathcal{A}}(\mathcal{B}(z)).
\end{eqnarray*}
As for the second formula, we set
$$\widetilde{L}_{\mathcal{A}}=\frac{\overline{c}}{\overline{c}\bsb{\xi}+\overline{d}},$$ where we recall that
$\bsb{\xi}=\frac{\bsb{s}_1}{\bsb{s}_2}$.
By using (\ref{eq:iotazeta}) and the obvious identity $\det(\mathcal{B})=\det(\ol{\mathcal{B}})$, we compute in a similar way:
\begin{eqnarray*}
\lefteqn{\bsb{J}_{\mathcal{B}}(z)^2(\widetilde{L}_{\mathcal{C}}(z)-\widetilde{L}_{\mathcal{B}}(z))=}\\
&=&(\ol{\gamma} z+\ol{\delta})^2\left(\frac{\ol{x}}{\ol{x}\bsb{\xi}+\ol{y}}-\frac{\ol{\gamma}}{\ol{\gamma}\bsb{\xi}+\ol{\delta}}\right)\\
&=&\det(\ol{\mathcal{B}})\frac{\ol{c}}{\ol{c}\frac{\ol{\alpha}\bsb{\xi}+\ol{\beta}}{\ol{\gamma}\bsb{\xi}+\ol{\delta}}+\ol{d}}\\
&=&\det(\mathcal{B})\widetilde{L}_{\mathcal{A}}(\mathcal{B}(z)).
\end{eqnarray*}
Hence,
$$\widetilde{L}_{\mathcal{A}}(\mathcal{B}(z))=\det(\mathcal{B})^{-1}\boldsymbol{J}_{\mathcal{B}}(z)^2(\widetilde{L}_{\mathcal{C}}(z)-\widetilde{L}_{\mathcal{B}}(z)).$$
But
$$\widetilde{L}_{\mathcal{A}}(z)=\bfs_2(z)\boldsymbol{L}_{\mathcal{A}}(z),$$
so that
\begin{eqnarray*}
\widetilde{L}_{\mathcal{A}}(\mathcal{B}(z))&=&\bfs_2(\mathcal{B}(z))\boldsymbol{L}_{\mathcal{A}}(\mathcal{B}(z))\\
&=&(\overline{\gamma}\bfs_1(z)+\overline{\delta}\bfs_2(z))\boldsymbol{L}_{\mathcal{A}}(\mathcal{B}(z))\\
&=&\bfs_2(z)J_{\mathcal{B}}(z)^{-1}\boldsymbol{J}_{\mathcal{B}}(z)\boldsymbol{L}_{\mathcal{A}}(\mathcal{B}(z)),
\end{eqnarray*}
where $\bfs_1,\bfs_2$ are considered as functions $\Omega\rightarrow\TT_{>0}$,
from which we deduce the expected identity.\CVD

\subsection{Almost $A$-quasi-modular forms.}

We recall that for all $z\in\Omega$ and $\gamma\in\Gamma$, we have $J_\gamma,\boldsymbol{J}_\gamma,L_\gamma,\boldsymbol{L}_\gamma\in \TT_{>0}$.

Let $r$ be a positive real number and $\bsb{f}:\Omega\rightarrow\TT_{<r}$ a map. 
We will say that $\bsb{f}$ is {\em regular} if the following properties hold.

\begin{enumerate}
\item There exists $\varepsilon>0$ such that, for all $t_0\in C$, $|t_0|<\varepsilon$,
the map $z\mapsto \bsb{f}(z,t_0)$ is
holomorphic on $\Omega$.

\item
For all $a\in A$, $\bsb{f}(z+a)=\bsb{f}(z)$. Moreover, there exists $c>0$ such that for all $z\in\Omega$ with
$|u(z)|<c$ and $t$ with $|t|<c$, there is a convergent expansion
\begin{equation*}\label{uexpansion}
\bsb{f}(z,t)=\sum_{n,m\geq 0}c_{n,m}t^nu^m,
\end{equation*}
where $c_{n,m}\in C$.
\end{enumerate}

\begin{Definition}[Almost $A$-quasi-modular forms]\label{def_almost}
{\em Let $\bsb{f}$ be a regular function $\Omega\rightarrow\TT_{<r}$, for $r$ a positive real number. 
We say that $\bsb{f}$ is an {\em almost-$A$-quasi-modular form of weight $(\mu,\nu)$,  
type $m$ and depth $\leq l$} if there exist regular functions $\bsb{f}_{i,j}:\Omega\rightarrow\TT_{<r}$,
$0\le i+j\le l$, such that for all $\gamma\in\Gamma$
and $z\in\Omega$ the following functional equation holds in $\TT_{>0}$:
\begin{equation}\label{functional:equations}
\bsb{f}(\gamma(z),t)=\det(\gamma)^{-m}J_\gamma^\mu\bsb{J}_\gamma^\nu\left(\sum_{i+j\leq l}\bsb{f}_{i,j}L_\gamma^i\bsb{L}_\gamma^j\right).
\end{equation}}
\end{Definition}

The {\em radius of convergence} $\rho(\bsb{f})$ of an almost $A$-quasi-modular form $\bsb{f}:\Omega\rightarrow\TT_{>0}$ is the supremum of the set of the real numbers 
$r$ such that the maps $\bsb{f},\bsb{f}_{i,j}$ appearing in (\ref{functional:equations}) simultaneously are well defined maps $\Omega\rightarrow\TT_{<r}$.

We will say that $\mu=\mu(\bsb{f}),\nu=\nu(\bsb{f}),m=m(\bsb{f})$
are respectively the {\em first weight}, the {\em second weight} and the {\em type} of $\bsb{f}$.
\subsubsection{Some remarks.}
It is obvious that in (\ref{functional:equations}),
$\bsb{f}=\bsb{f}_{0,0}$ (use $\gamma=$ identity matrix).

If $\lambda\in\TT_{>0}$, then the map $z\mapsto\lambda$ trivially is an almost $A$-quasi-modular form of weight $(0,0)$, type $0$, depth $\leq 0$. The radius
$\rho(\lambda)$ is then just the radius of convergence of the series $\lambda$. 

Examples of almost $A$-quasi-modular forms are Drinfeld quasi-modular forms.
To any Drinfeld quasi-modular form of weight $w$, type $m$, depth $\leq l$ is associated an almost $A$-quasi-modular form of weight $(w,0)$, type $m$, 
depth $\leq l$ whose radius is infinite. 

The $\TT_{>0}$-algebra $\TT_{>0}[g,h]$ is graded by 
the couples $(w,m)\in\ZZ\times\ZZ/(q-1)\ZZ$ of weights and types, and the isobaric elements are all almost $A$-quasi-modular forms
with the second weight $0$.

The function $\bfs_2$ is, by Lemmas \ref{lemmavectorial} and \ref{usefulford}, an almost $A$-quasi-modular form of weight $(-1,1)$, depth $0$, type $0$.
The radius is $q$, by the results of Section \ref{edm}.

If $\bsb{f}$ is an almost $A$-quasi-modular form of weight $(\mu,\nu)$, type $m$, depth $\leq l$ and radius of convergence $>q$, then
$\varepsilon(\bsb{f}):=\bsb{f}|_{t=\theta}$ is a well defined holomorphic function $\Omega\rightarrow C$. 
It results from (\ref{Jbold}) and (\ref{Lbold}) that $\varepsilon(\bsb{f})$ is a Drinfeld quasi-modular form of weight $\mu+\nu$, type $m$ and depth $\leq l$.

The function $\bsb{f}:=\bfs_2$ is not well defined at $t=\theta$ because its radius of convergence is $q$, and we know from (\ref{uexpansions2}) that there is divergence at $\theta$. However, the function
$\bsb{f}:=(t-\theta)\bfs_2$, which is an almost $A$-quasi-modular of same weight, type and depth as $\bfs_2$, has convergence radius $q^q$. Therefore,
$\varepsilon(\bsb{f})$ is well defined, and is the constant function $-1$ by the results of Subsection \ref{functionsofthevariable}.
From (\ref{functionalequationstwist}) we see that the function $\bsb{s}_2^{(1)}$ is not an almost $A$-quasi-modular form. The non-zero function
$\varepsilon(\bsb{s}_2^{(1)})$ is well defined and we have already mentioned the results of Gekeler in \cite{gekeler:compositio} that allow to compute it.

Let us write $\phi=\varepsilon(\bsb{E})$, which corresponds to a well defined series of $u C[[u^{q-1}]]$ by  (\ref{eq:u-expansionofE}). 
We obtain, by using (\ref{tmodularityE}), (\ref{Jbold}) and (\ref{Lbold}) with $k=0$, that
$$\phi(\gamma(z))=\det(\gamma)^{-1}(cz+d)^{2}\left(\phi(z)-\widetilde{\pi}^{-1}\frac{c}{cz+d}\right).$$ This is the collection of functional equations
of the Drinfeld quasi-modular form $E$ (\ref{formE}), whose $u$-expansion begins with the term $u$. 
Applying \cite[Theorem 1]{BP} we obtain:

\begin{Lemme}\label{EEgivesE}
We have, for all $z\in\Omega$:
$$\varepsilon(\bsb{E})=\bsb{E}(z,\theta)=E(z).$$
\end{Lemme}

It is easy to verify, as a confirmation of this result, that the first coefficients of the $u$-expansion of $\boldsymbol{E}$ given in (\ref{eq:u-expansionofE})
agree, substituting $t$ by $\theta$, with the $u$-expansion of $E$ that we know already after 
\cite[Corollary (10.5)]{Ge}:
\begin{equation*}
E=u(1+v^{(q-1)}+\cdots).
\end{equation*}

More generally, Propositions \ref{taudifferenceeq}, \ref{propositionautomorphyE} and \ref{propositionEuexpansion}
imply that for all $k\geq 0$, $\bsb{E}^{(k)}$ is an almost $A$-quasi-modular form of weight $(q^k,1)$ type $1$ and depth $\leq 1$ with convergence radius $\geq q^q$,
so that $\varepsilon(\bsb{E}^{(k)})$ is well defined, and is a Drinfeld quasi-modular form of weight $q^k+1$, type $1$ and depth $\leq 1$.

\subsubsection{Grading by the weights, filtering by the depths.}

For $\mu,\nu\in\ZZ,m\in\ZZ/(q-1)\ZZ$, $l\in\ZZ_{\geq 0}$, we denote by $\widetilde{\mathcal{M}}_{\mu,\nu,m}^{\leq l}$ the $\TT_{>0}$-module of
almost $A$-quasi-modular forms of weight $(\mu,\nu)$, type $m$ and depth $\leq l$. 
We have $$\widetilde{\mathcal{M}}^{\leq l}_{\mu,\nu,m}\widetilde{\mathcal{M}}^{\leq l'}_{\mu',\nu',m'}\subset\widetilde{\mathcal{M}}^{\leq l+l'}_{\mu+\mu',\nu+\nu',m+m'}.$$
We also denote by $\widetilde{\mathcal{M}}$ the $\TT_{>0}$-algebra
generated by all the almost $A$-quasi-modular forms. We prove below that this algebra is graded by the group $G=\ZZ^2\times \ZZ/(q-1)\ZZ$, filtered by the depths
(Proposition \ref{gradingfiltering}),
and contains five algebraically independent functions $E,g,h,\bsb{E},\bsb{F}$ (Proposition \ref{algebraicindependenceEghE1}).

Let $\mathcal{K}$ be any field extension of  $\FF_q(t,\theta)$. The key result of this section is the following elementary lemma.
\begin{Lemme}\label{zariski}
The subset $\Theta=\{(d,\overline{d}),d\in A\}\subset\A^2(\mathcal{K})$ is Zariski dense. 
\end{Lemme}
\noindent\emph{Proof.} 
Let us assume by contradiction that the lemma is false and let $\overline{\Theta}$ be the Zariski closure of $\Theta $.
Then, we can write 
$$\overline{\Theta}=\bigcup_{i\in\mathcal{I}} \Theta_i\cup\bigcup_{j\in\mathcal{J}}\widetilde{\Theta}_j,$$
where the $\Theta_i$'s are irreducible closed subsets of $\A^2(\mathcal{K})$ of dimension $1$, 
the $\widetilde{\Theta}_j$'s are isolated points of $\A^2(\mathcal{K})$, and $\mathcal{I},\mathcal{J}$ are finite sets.

From $\Theta = \Theta +(d,\ol{d})$ for all $d\in A$ we deduce $\ol{\Theta}=\ol{\Theta}+(d,\ol{d})$. The translations 
of $\A^2(\mathcal{K})$ by points such as $(d,\ol{d})$ being bijective, they induce permutations of the sets $\{\Theta_i\}$
and $\{\widetilde{\Theta}_j\}$, from which we easily deduce that $\mathcal{J}=\emptyset$.
Therefore, the ideal of polynomials $R\in\mathcal{K}[X,Y]$ such that $R(\Theta)\subset\{0\}$ is principal, generated by a
non-zero polynomial $P$.

Now, if $b\in A$, $m_b(\ol{\Theta})\subset\ol{\Theta}$, where $m_b(x,y):=(bx,\ol{b}y)$. Hence, $P(m_b(X,Y))\in(P)$ and
there exists $\kappa_b\in\mathcal{K}^\times$ such that 
$$P(bX,\ol{b}Y)=\kappa_bP(X,Y).$$
Let us write:
$$P(X,Y)=\sum_{\alpha,\beta}c_{\alpha,\beta}X^\alpha Y^\beta,$$ and choose $b\not\in\FF_q$.
If $c_{\alpha,\beta}\neq 0$, then $\kappa_b=b^{-\alpha}\ol{b}^{-\beta}$. If $P$ is not a monomial, we have, for $(\alpha,\beta)\neq(\alpha',\beta')$, $c_{\alpha,\beta},c_{\alpha',\beta'}\neq 0$, so that $b^{-\alpha}\ol{b}^{-\beta}=b^{-\alpha'}\ol{b}^{-\beta'}$,
yielding a contradiction, because $b\not\in\FF_q$.

If $P$ is a monomial, however, it cannot vanish at $(1,1)\in \Theta $; contradiction.\CVD

\begin{Lemme}\label{vanishingJJ}
Let us suppose that for elements $\psi_{\alpha,\beta}\in C((t))$ and for a certain element $z\in\Omega$ we have an identity:
\begin{equation}\label{eq:relationalphabeta}
\sum_{\alpha,\beta}\psi_{\alpha,\beta}J_\gamma^\alpha \boldsymbol{J}_\gamma^\beta=0,
\end{equation}
in $C((t))$, for all $\gamma=\begin{pmatrix}a & b \\ 1 & d\end{pmatrix}\in\Gamma$ with determinant $1$, the sum being finite.
Then, $\psi_{\alpha,\beta}=0$ for all $\alpha,\beta$.
\end{Lemme}
\noindent\emph{Proof.} 
Let us suppose by contradiction the existence of a non-trivial relation (\ref{eq:relationalphabeta}).
We have, with the hypothesis on $\gamma$, $J_\gamma=z+d,\boldsymbol{J}_\gamma=
\bsb{\xi}+\overline{d}\in C((t))$, so that 
the relation of the lemma implies the existence of a relation:
$$\sum_{\alpha,\beta}\ell_{\alpha,\beta}d^\alpha\overline{d}^\beta=0,\quad d\in A,$$
with $\ell_{\alpha,\beta}\in\mathcal{K}=C((t))$ not all zero, and all, but finitely many, vanishing. Lemma \ref{zariski}
yields a contradiction.\CVD

Another useful lemma is the following. The proof is again a simple application of Lemma \ref{zariski} and
will be left to the reader.

\begin{Lemme}\label{lemmaLABC}
If the finite collection of functions $f_{i,j}:\Omega\rightarrow \TT_{>0}$ is such that for all $z\in\Omega$ and for all $\gamma\in\Gamma$,
$$\sum_{i,j}f_{i,j}(z)L_\gamma^i\boldsymbol{L}_\gamma^j=0,$$
then the functions $f_{i,j}$ are all identically zero.
\end{Lemme}

\begin{Lemme}\label{lemme_types}
Let $\bsb{f}$ be an almost $A$-quasi-modular form of type $m$ with $0\leq m<q-1$. Then, with $v=u^{q-1}$,
$$\bsb{f}(z)=u^m\sum_{i\geq 0}c_i(t)v^{i}.$$
\end{Lemme}
\noindent\emph{Proof.} It follows the same ideas of the remark on p. 23
of \cite{Ge1}. Let us consider $\gamma=\begin{pmatrix}\lambda & 0 \\ 0 & 1\end{pmatrix}\in\Gamma$ with $\lambda\in\FF_q^\times$. We have 
$\gamma(z)=\lambda z$, $\det(\gamma)=\lambda$, $J_\gamma=\bsb{J}_\gamma=1$, $L_\gamma=\bsb{L}_\gamma=0$, so that
$\bsb{f}(\lambda z)=\lambda^{-m}\bsb{f}(z)$, for all $z\in\Omega$.
Now, if $\bsb{f}=\sum_ic_i(t)u^i$, since $e_\carlitz$ is $\FF_q$-linear, we get $u(\lambda z)=\lambda^{-1}u(z)$ and 
if $c_i\neq0$, then $i\equiv m\pmod{q-1}$.\CVD

\begin{Proposition}\label{gradingfiltering}
The $\TT_{>0}$-algebra generated by the almost $A$-quasi-modular forms is graded by weights and types, 
hence by the group $G=\ZZ^2\times\ZZ/(q-1)\ZZ$, and filtered by the depths:
$$\widetilde{\mathcal{M}}=\bigoplus_{(\mu,\nu,m)\in G}\;\bigcup_{l=0}^\infty\widetilde{\mathcal{M}}_{\mu,\nu,m}^{\leq l}.$$
\end{Proposition}

\noindent\emph{Proof.} We begin by proving the property concerning the grading by the group $\ZZ^2\times\ZZ/(q-1)\ZZ$. 
Let us consider distinct triples $(\mu_i,\nu_i,m_i)\in\ZZ^2\times\ZZ/(q-1)\ZZ$, $i=1,\ldots,s$, non-negative integers $l_1,\ldots,l_s$  and non-zero elements $\bsb{f}_i\in\widetilde{\mathcal{M}}^{\leq l_i}_{\mu_i,\nu_i,m_i}$. Then, we claim that $\sum_{i=1}^s\bsb{f}_i\neq 0$.
To see this, we assume by contradiction that for some forms $\bsb{f}_i$ as in the proposition, we have the identity in $\TT_{>0}$:
\begin{equation}\label{hypothesis1}
\sum_{i=1}^s\bsb{f}_i=0.
\end{equation}
Recalling Definition \ref{def_almost} (identity (\ref{functional:equations})), we have, 
for all $i=1,\ldots,s$, $\gamma=\begin{pmatrix}a & b \\ c & d\end{pmatrix}\in\Gamma$, $z\in\Omega$:
$$\bsb{f}_i(\gamma(z),t)=\det(\gamma)^{-m_i}J_\gamma^{\mu_i}\boldsymbol{J}_\gamma^{\nu_i}\sum_{j+k\leq l}\bsb{f}_{i,j,k}(z,t)L_\gamma^j\boldsymbol{L}_\gamma^k,$$
for certain functions $\bsb{f}_{i,j,k}:\Omega\rightarrow\TT_{>0}$.

Let us suppose first that $\gamma$ is of the form $\begin{pmatrix} a & b\\ 1 & d\end{pmatrix}$ with 
$ad-b=1$. We recall that $\bsb{s}_2(z)^{-1}\in\TT_{>0}^\times$ for all $z$. Therefore, for all $z\in\Omega$, (\ref{hypothesis1}) becomes the identity of formal series in $\TT_{>0}$:
\begin{equation}\label{hypothesis2}
\sum_{i=1}^s\sum_{j+k\leq l_i}\bsb{f}_{i,j,k}\boldsymbol{s}_2^{-k}(z+d)^{\mu_i-j}(\bsb{\xi}+\overline{d})^{\nu_i-k}=0.
\end{equation}
By Lemma \ref{vanishingJJ}, (\ref{hypothesis2}) is equivalent to the relations:
\begin{equation}\label{alphabeta}
\sum_{i,j,k}\phi_{i,j,k}=0,\quad\text{ for all }(\alpha,\beta)\in\ZZ^2
\end{equation}
where $\phi_{i,j,k}:=\bsb{f}_{i,j,k}\boldsymbol{s}_2^{-k}$ and the sum runs over the triples
$(i,j,k)$ with $i\in\{1,\dots,s\}$ and $j,k$ such that $\mu_i-j=\alpha$ and $\nu_i-k=\beta$, 
with obvious vanishing conventions on some of the $\phi_{i,j,k}$'s.

Let $\mu$ be the maximum value of the $\mu_i$'s, and let us look at the relations (\ref{alphabeta}) for 
$\alpha=\mu$. Since for all $\mu_i<\mu$ we get $\alpha=\mu>\mu_i-j$ for all $j\geq 0$,
for such a choice of $\alpha$ we get:
\begin{equation}\label{alphabeta2}
\sum_{i,k}\phi_{i,0,k}=0,\quad\text{ for all }\beta\in\ZZ,
\end{equation}
where the sum is over the couples $(i,j)$ with $i$ such that $\mu_i=\mu$ and $\nu_i-k=\beta$.
Now, let $\mathcal{E}$ be the set of indices $i$ such that $\mu_i=\mu$ and write 
$\nu$ for the maximum of the $\nu_i$ with $i\in\mathcal{E}$. If $j$ is such that $\mu_j=\mu$,
and if $\nu\neq\nu_j$, then for all $k\geq 0$, $\nu>\nu_j-k$, so that for $\beta=\nu$, 
(\ref{alphabeta2}) becomes 
$$\sum_{i}\phi_{i,0,0}=0,$$ where the sum runs this time over the $i$'s such that
$(\mu_i,\nu_i)=(\mu,\nu)$.
But $\phi_{i,0,0}=\bsb{f}_{i,0,0}=\bsb{f}_i$ for $i=1,\ldots,s$. Since the types of the $\bsb{f}_i$'s with same
weights are distinct by hypothesis, Lemma \ref{lemme_types} implies that 
for all $i$ such that $(\mu_i,\nu_i)=(\mu,\nu)$, $\bsb{f}_i=0$. This contradicts our initial assumptions and proves our initial claim.
Combining with Lemma \ref{lemmaLABC}, we end the proof of the proposition.\CVD

\begin{Proposition}\label{algebraicindependenceEghE1}
The functions $$E,g,h,\bfs_2,\bfs_2^{(1)}:\Omega\rightarrow \TT_{>0}$$
are algebraically independent over the fraction field of $\TT_{>0}$.
\end{Proposition}

\noindent\emph{Proof.} Assume by contradiction that the statement of the proposition is false. Since $E,g,h,\bsb{s}_2,\bsb{s}_2^{(1)}\in\widetilde{\mathcal{M}}$ are
almost $A$-quasi-modular forms, by Proposition \ref{gradingfiltering}, 
there exist $(\mu,\nu),m\in\ZZ$, and a non-trivial relation (where the sum is finite):
$$\sum_{i,j\geq 0}P_{i,j}E^i\bfs_2^{(1)}{}^j=0,$$
with $P_{i,j}\in\TT_{>0}[g,h,\bsb{s}_2]\cap\widetilde{\mathcal{M}}^{\leq l}_{\mu-2i+qj,\nu-j,m-i}$ (for some $l\geq 0$). By Proposition \ref {gradingfiltering},
any vector space of almost $A$-quasi-modular forms of given weight and depth is filtered by the depths. Comparing with the functional equations (\ref{functionalequationstwist}) and \cite[Functional equation (11)]{BP}, and applying Lemma \ref{lemmaLABC}, we see that
all the forms $P_{i,j}$ vanish. There are three integers $\alpha,m,n$ and a non trivial polynomial relation $P$ among $g,h,\bfs_2$, with coefficients in $\TT_{>0}$:
$$\sum_{s=0}^nQ_s\bfs_2^s=0,$$
where $Q_s\in \TT_{>0}[g,h]\cap\widetilde{\mathcal{M}}^{\leq l}_{\alpha+s,0,m}$ ($s=0,\ldots,n$), and for some $s$, $Q_s$ is non-zero.
Since $\nu(Q_s)=0$ for all $s$ such that $Q_s\neq0$ and $\nu(\bsb{s}_2)=1$,
The polynomial $P$, evaluated at the functions $E,g,h,\bsb{s}_2,\bsb{s}_2^{(1)}$ is equal to $Q\bsb{s}_2^s$ for $Q\in \TT_{>0}[g,h]\setminus\{0\}$ and $s\in\ZZ$, quantity that cannot vanish because $g,h$ are algebraically independent over $\TT_{<q^q}$: contradiction.\CVD

\section{Estimating the multiplicity\label{estimating}}

We prove Theorem \ref{secondtheorem} in this section.
\subsection{Preliminaries}

Let us denote by $\MM^\dag$ the $\TT_{<q^q}$-algebra $\TT_{<q^q}[g,h,\boldsymbol{E},\boldsymbol{F}]$, where $\boldsymbol{F}:=\boldsymbol{E}^{(1)}$; its dimension is $4$, according to Proposition \ref{algebraicindependenceEghE1} and Proposition \ref{propositionautomorphyE}. By Proposition \ref{gradingfiltering}, this algebra is graded by the group  
$\ZZ^2\times\ZZ/(q-1)\ZZ$:
$$\MM^\dag=\bigoplus_{(\mu,\nu),m}\MM^\dag_{\mu,\nu,m},$$ where $\MM^\dag_{\mu,\nu,m}=\widetilde{\mathcal{M}}_{\mu,\nu,m}\cap\MM^\dag$.

The operator $\tau$ acts on $\MM^\dag$ by Proposition \ref{taudifferenceeq}. More precisely, 
we have the homomorphism of $\FF_q[t]$-modules $$\tau:\MM^\dag_{\mu,\nu,m}\rightarrow\MM^\dag_{q\mu,\nu,m}.$$

Let us write $\bsb{h}=\widetilde{\pi}hs_\carlitz^{-1} \bsb{s}_2=h\bsb{d}$.

\begin{Lemme}\label{h_1}
The formula $\boldsymbol{h}=(t-\theta^q)\boldsymbol{F}-g\boldsymbol{E}$ holds, so that $\bsb{h}\in\MM^\dag_{q,1,1}$ and $\MM^\dag=\TT_{<q^q}[g,h,\bsb{E},\bsb{h}]$.
\end{Lemme}
\noindent\emph{Proof.} From the definition of $\boldsymbol{E}$ and (\ref{equexpansion}), we find:
\begin{eqnarray*}
\lefteqn{(t-\theta^q)\boldsymbol{F}-g\boldsymbol{E}=}\\&=&-(t-\theta^q)h^q\bsb{d}^{(2)}+gh\bsb{d}^{(1)}\\
&=&(-h)^q(-h^{q-1})^{-1}(\bsb{d}-g\bsb{d}^{(1)})+gh\bsb{d}^{(1)}\\
&=&h\bsb{d}=\bsb{h}.
\end{eqnarray*}
This makes it clear that $\bsb{h}$ belongs to $\MM^\dag_{q,1,1}$ and that $\MM^\dag=\TT_{<q^q}[g,h,\bsb{E},\bsb{h}]$.\CVD

We denote by $\varepsilon_{\mu,\nu,m}$ or again $\varepsilon$ the map
which sends an almost $A$-quasi-modular form $\bsb{f}$ of weight $(\mu,\nu)$, type $m$, with radius $>q$ to the 
Drinfeld quasi-modular form $\varepsilon(\bsb{f})$ of weight $\mu+\nu$, type $m$. This map is clearly a $C$-algebra homomorphism.

\begin{Lemme}\label{defh=h1}
We have $\varepsilon(\bsb{h})=h$.
\end{Lemme}
\noindent\emph{Proof.} This follows from the limit $\lim_{t\rightarrow\theta}s_\carlitz^{-1}\bsb{s}_2=\widetilde{\pi}^{-1}$
and the definition of $\bsb{d}$.\CVD

More generally, we have the following result.

\begin{Proposition}\label{almost_bijective} For all $(\mu,\nu),m$, 
the map $$\varepsilon:\MM^\dag_{\mu,\nu,m}\rightarrow\widetilde{M}^{\leq \nu}_{\mu+\nu,m}$$ is well defined and the inverse image of $0$
is the $\TT_{<q^q}$-module $(t-\theta)\MM^\dag_{\mu,\nu,m}$.
\end{Proposition}

\noindent\emph{Proof.} 
Let $\bsb{f}$ be an element of $\MM^\dag_{\mu,\nu,m}$. Then, by Lemma \ref{h_1},
$$\bsb{f}=\sum_{i=0}^\nu\phi_i\bsb{h}^{\nu-i}\bsb{E}^i,$$ where $\phi_i\in M_{\mu-\nu q+i(q-1),m-\nu}\otimes_C\TT_{<q^q}$.
Since $\lim_{t\rightarrow\theta}s_\carlitz^{-1}\bsb{s}_2=\widetilde{\pi}^{-1}$, we have $\varepsilon(\bsb{h})=h$ by Lemma \ref{defh=h1}. Moreover, by Lemma \ref{EEgivesE}, $\varepsilon(\bsb{E})=E$, and
\begin{eqnarray*}
\varepsilon(\bsb{f})&=&\sum_{i=0}^\nu\varepsilon(\phi_i)h^{\nu-i}E^i,
\end{eqnarray*}
so that $\varepsilon(\bsb{f})=0$ if and only if $\varepsilon(\phi_i)=0$ for all $i$. But for all $i$, $\phi_i$ is a polynomial in $g,h$ with coefficients in $\TT_{<q^q}$.
If $\varepsilon(\phi_i)=0$, then $\phi_i$ is a linear combination $\sum_{a,b}c_{a,b}g^ah^b$ with $c_{a,b}\in\TT_{<q^q}$ such that $c_{a,b}(\theta)=0$.
Since $\TT_{<q^q}\subset\TT$, it is a principal ideal domain and the last condition is equivalent to $\phi_i\in(t-\theta)(M\otimes_C\TT_{<q^q})$.
Hence, $\varepsilon(\bsb{f})=0$ if and only if, for all $i$, $\phi_i\in(t-\theta)(M\otimes_C\TT_{<q^q})$. The proposition follows.\CVD

\subsection{Multiplicity estimate in $\MM^\dag$\label{multiMMdag}}

By Proposition \ref{propositionEuexpansion}, $\bsb{E}=u+\cdots\in u\FF_q[t,\theta][[u^{q-1}]]$. 
Hence,
$$\bsb{E}^{(k)}=u^{q^k}+\cdots\in u^{q^k}\FF_q[t,\theta][[u^{(q-1)q^k}]],\quad k\geq 0,$$
and there is an embedding $\MM^\dag\rightarrow \TT_{<q^q}[[u]]$. It will be sometimes useful to fix an embedding of
$\TT_{<q^q}$ in $\mathcal{K}$, an algebraic closure of $C((t))$; we will then often consider elements of $\MM^\dag$ as 
formal series if $\mathcal{K}[[u]]$ (especially in this subsection). Anderson's 
operator $\tau:C((t))\rightarrow C((t))$ extends in a natural way to an $\FF_q(t)$-linear operator $\tau:\mathcal{K}\rightarrow\mathcal{K}$
(we will keep using the notation $\tau^k f=f^{(k)}$).
If $f=\sum_{n\geq n_0}c_n(t)u^{n}$ is a formal series of $\mathcal{K}[[u]]$, then, Anderson's operator further extends as follows:
\begin{equation}\label{ftauk}
f^{(k)}=\sum_{n\geq n_0}c_n^{(k)}(t)u^{q^kn},\quad k\in\ZZ.
\end{equation}
Let $f=\sum_{n\geq n_0}c_n(t)u^{n}$ be in $\mathcal{K}[[u]]$, with $c_{n_0}\neq 0$. We write $\nu_\infty(f):=n_0$. We also set $\nu_\infty(0):=\infty$. 
Obviously, $\nu_{\infty}(f^{(k)})=q^k\nu_{\infty}(f)$ for all $k\geq 0$.
We recall that $\nu_\infty(g)=0,\nu_\infty(h)=\nu_\infty(\bsb{E})=1$ and $\nu_\infty(\bsb{F})=q$. Since $\nu_\infty(\bsb{s}_2)=0$, we also get $\nu_\infty(\bsb{h})=1$.
In the following, we will write $\MM^\dag_{\mu,\nu,m}(\mathcal{K})=\MM^\dag_{\mu,\nu,m}\otimes_{\TT_{<q^q}}\mathcal{K}$ and $M_{w,m}(\mathcal{K})=M_{w,m}\otimes_C\mathcal{K}$. It is evident that the $\mathcal{K}$-algebra $\MM^\dag(\mathcal{K})=\sum_{\mu,\nu,m}\MM^\dag_{\mu,\nu,m}(\mathcal{K})$ is again graded by the group $\ZZ^2\times\ZZ/(q-1)\ZZ$;
similarly for the algebra $M(\mathcal{K})=\sum_{w,m}M_{w,m}(\mathcal{K})$.

We begin with a rather elementary estimate, for $\bsb{f}\in\MM^\dag$ of weight $(\mu,0)$.

\begin{Lemme}\label{lemmeriemannroch}
If $\bsb{f}\in\MM^\dag_{\mu,0,m}(\mathcal{K})$ is non-zero, then $\nu_\infty(\bsb{f})\leq \frac{\mu}{q+1}$.
\end{Lemme}
\noindent\emph{Proof.} 
A weight inspection shows that $\MM^\dag_{\mu,0,m}(\mathcal{K})=\mathcal{K}[g,h]_{\mu,m}$.
We can write $\bsb{f}=h^{\nu_\infty(\bsb{f})}\bsb{b}$, with $\bsb{b}\in\mathcal{K}[g,h]$ and $h$ not dividing $\bsb{b}$. Therefore, $\nu_\infty(\bsb{f})\leq \frac{\mu}{q+1}$.\CVD

In the next proposition, we study the case of $\bsb{f}$ of weight $(\mu,\nu)$ with $\nu>0$.
\begin{Proposition}\label{Mdag}
Let $\bsb{f}$ be a non-zero element of $\MM^\dag_{\mu,\nu,m}(\mathcal{K})$ with $\nu\neq 0$. Then, 
$$\nu_\infty(\bsb{f})\leq \mu \nu.$$
\end{Proposition}
 It is not difficult to show that the statement of this proposition cannot be improved (this can be checked with the functions
 $\bsb{E}^{(k)}$ in mind).

Before proving the proposition, we need to state and prove a lemma.

\begin{Lemme}\label{resultantffprimo}
Let $\bsb{f}\in\MM^\dag_{\mu,\nu,m}(\mathcal{K})$, $\bsb{f}'\in\MM^\dag_{\mu',\nu',m'}(\mathcal{K})$. By Lemma \ref{h_1}, $\bsb{f},\bsb{f}'$ can be written in an unique way as polynomials in $\mathcal{K}[g,h,\bsb{E},\bsb{h}]$. Let $l,l'$ be the degrees in $\bsb{E}$ of $\bsb{f},\bsb{f}'$  respectively. Then (Resultant), 
$$\bsb{\phi}:=\mathbf{Res}_{\bsb{E}}(\bsb{f},\bsb{f}')=\bsb{h}^{\nu l'+\nu' l-ll'}\phi_0,$$ where
$\phi_0\in M_{w^*,m^*}(\mathcal{K})$, with
$$w^*=\mu l'+\mu' l-ll'-q(\nu l'+\nu' l-ll'),\quad m^*:=m l'+m' l-(\nu l'+\nu' l).$$
\end{Lemme}

\noindent\emph{Proof.} With an application of an obvious variant of \cite[Lemme 6.1]{PF} (\footnote{The first formula after the statemement of the above cited lemma, mistakenly typed, must be replaced 
with \[p(R)=p(F)\deg_{X_0}(G)+p(G)\deg_{X_0}(F)-p(X_0)\deg_{X_0}(F)\deg_{X_0}(G).\]}) we see that
$$\bsb{\phi}\in\MM^\dag_{\mu l'+\mu' l-ll',\nu l'+\nu' l-ll',m l'+m' l-ll'}(\mathcal{K}).$$
At the same time, $\bsb{\phi}\in\mathcal{K}[g,h,\bsb{h}]$. Since
$\nu(g)=\nu(h)=0$ and $\nu(\bsb{h})=1$, we have $\phi_0:=\bsb{\phi}/\bsb{h}^{\nu l'+\nu' l-ll'}\in M(\mathcal{K})$.
The computation of the weight and type of $\phi_0$ is obvious, knowing that $\mu(\bsb{h})=q$.\CVD

\noindent\emph{Proof of Proposition \ref{Mdag}.} 
Let $\bsb{f}$ be in $\MM^\dag_{\mu,\nu,m}(\mathcal{K})$, with $\nu>0$. Assume first that $\bsb{f}$, as a polynomial in $g,h,\bsb{E},\bsb{h}$, is irreducible. 
If $\bsb{f}$ belongs to $\mathcal{K}[g,h,\bsb{h}]$ then $\bsb{f}=\phi\bsb{h}^\nu$ with $\phi\in M_{\mu-q\nu,0,m-\nu}(\mathcal{K})$ and 
\begin{eqnarray*}
\nu_\infty(\bsb{f})&\leq&\nu_{\infty}(\phi)+\nu\nu_{\infty}(\bsb{h})\\
&\leq &\frac{\mu-q\nu}{q+1}+\nu\\
&\leq &\frac{\mu+\nu}{q+1}\\
&\leq &\mu\nu.
\end{eqnarray*}
We now suppose that $\bsb{f}\not\in\mathcal{K}[g,h,\bsb{h}]$; there are two cases left.

\medskip

\noindent\emph{Case (i).} We suppose that $\bsb{f}$ divides $\bsb{f}^{(1)}\in\MM^\dag_{q\mu,\nu,m}(\mathcal{K})$ as a polynomial in $g,h,\bsb{E},\bsb{h}$. For weight reasons, $\bsb{f}^{(1)}=a\bsb{f}$ with $a\in M_{\mu(q-1),0}(\mathcal{K})$ and $a\neq 0$. We also have $\nu_{\infty}(\bsb{f}^{(1)})=q\nu_\infty(\bsb{f})$ by (\ref{ftauk}), so that, by Lemma \ref{lemmeriemannroch}, $(q-1)\nu_\infty(\bsb{f})=\nu_\infty(a)\leq (q-1)(q+1)^{-1}\mu$.
Hence, in this case, we get the stronger inequality (\footnote{It can be proved that $\bsb{f}$ is, in this case, a modular form multiplied by an element of $\mathcal{K}$, but we do not need this information here.}) $$\nu_\infty(\bsb{f})\leq \frac{\mu}{q+1}.$$

\medskip

\noindent\emph{Case (ii).} In this case, $\bsb{f}$ and $\bsb{f}^{(1)}$ are
 relatively prime. Since $\bsb{f}$ is irreducible, $\deg_{\bsb{E}}(\bsb{f})=l=\nu>0$, so that $\bsb{f},\bsb{f}^{(1)}$ depend on $\bsb{E}$, and their resultant $\bsb{\phi}$ with respect to $\bsb{E}$ is non-zero. We apply Lemma \ref{resultantffprimo} with $\bsb{f}'=
\bsb{f}^{(1)}$, finding $$\bsb{\phi}=\bsb{h}^{\nu^2}\phi_0,$$
with $\phi_0\in M_{(q+1)\nu(\mu-\nu),m^*}(\mathcal{K})$, for a certain $m^*$ that can be computed with Lemma \ref{resultantffprimo}. By Lemma \ref{lemmeriemannroch} again, $\nu_\infty(\phi_0)\leq \nu(\mu-\nu)$.
Since $\nu_\infty(\bsb{h})=1$, $\nu_\infty(\bsb{\phi})\leq \nu(\mu-\nu)+\nu^2=\mu\nu$. Now, the number $\nu_\infty(\bsb{\phi})$
is an upper bound for $\nu_\infty(\bsb{f})$ by B\'ezout identity for the resultant.

\medskip

We have proved the proposition if $\bsb{f}\in\MM^\dag_{\mu,\nu,m}(\mathcal{K})$ is irreducible. If $\bsb{f}$ is not irreducible, we can write 
$\bsb{f}=\prod_{i=0}^r\bsb{f}_i$ with $\bsb{f}_0\in\MM^\dag_{\mu_0,0,m_0}(\mathcal{K})$, $\bsb{f}_i\in\MM^\dag_{\mu_i,\nu_i,m_i}(\mathcal{K})$ irreducible for all $i>0$ with $\nu_i>0$,
and $\sum_i\mu_i=\mu,\sum_i\nu_i=\nu,\sum_im_i\equiv m\pmod{q-1}$. Since $\nu_\infty(\bsb{f})=\sum_i\nu_\infty(\bsb{f}_i)$, we get, applying Lemma \ref{lemmeriemannroch},
$$\nu_\infty(\bsb{f})\leq\frac{\mu_0}{q+1}+\sum_{i>0}\mu_i\nu_i\leq\mu\nu.$$\CVD
\subsection{Reduced forms\label{reducedforms}}

Let $\bsb{f}$ be in $\MM^\dag$. Since $\varepsilon(\bsb{f})\in\widetilde{M}\subset C[[u]]$, it is legitimate 
to compare the quantities $\nu_\infty(\bsb{f})$ and $\nu_\infty(\varepsilon(\bsb{f}))$. We have the inequality:
\begin{equation}
\nu_\infty(\bsb{f})\leq\nu_\infty(\varepsilon(\bsb{f})),\label{nunu}
\end{equation}
but the equality is not guaranteed in general, because the leading term of the $u$-expansion of $\bsb{f}$ can vanish at $t=\theta$.

\begin{Definition}{\em A function $\bsb{f}$ in $\MM^\dag$ is {\em reduced} if $\nu_\infty(\bsb{f})=\nu_\infty(\varepsilon(\bsb{f}))$, that is, if the leading coefficient of the $u$-expansion of $\bsb{f}$ does not vanish at $t=\theta$.}
\end{Definition}


The next lemma provides a tool to construct reduced almost $A$-quasi-modular forms, useful in the sequel.

\begin{Lemme}\label{reduced} Let $\bsb{f}\in\MM^\dag_{\mu,\nu,m}$ be such that $\bsb{f}=\sum_{n\geq n_0}b_nu^n$, with $b_n\in\FF_q[t,\theta]$ for all $n$ and
$b_{n_0}\neq 0$. Then, for all $k> \log_q(\deg_tb_{n_0})$, the function $\bsb{f}^{(k)}$ is reduced.\end{Lemme}

\noindent\emph{Proof.} We have $b_{n_0}^{(k)}(\theta)=b_{n_0}(\theta^{q^{-k}})^{q^k}=0$ if and only if $t-\theta^{1/q^k}$ divides the polynomial $b_{n_0}(t)$
in $K^{\text{{alg.}}}[t]$. This polynomial having coefficients in $K$, we have $b_{n_0}^{(k)}(\theta)=0$ if and only if 
the irreducible polynomial $t^{q^k}-\theta$ divides $b_{n_0}(t)$. However, this is impossible if $k> \log_q(\deg_tb_{n_0})$.\CVD

\subsection{Construction of the auxiliary forms.}

We recall the $u$-expansion of $\bsb{E}$ whose existence is proved in Proposition \ref{propositionEuexpansion}:
$$\bsb{E}=u\sum_{i\geq 0}c_i(t)v^{i},$$ where $c_0=1$, $c_i\in\FF_q[t,\theta]$ for all $i>0$ and $v=u^{q-1}$.

\begin{Proposition}\label{propositionexpansions} The following properties hold.

\begin{description}
\item[(i)] Let $\alpha,\beta,\gamma,\delta$ be non-negative integers and let us write $\bsb{f}=g^\alpha h^\beta \bsb{E}^\gamma \bsb{F}^\delta\in\MM^\dag_{\mu,\nu,m}$,
with $\mu=\alpha(q-1)+\beta(q+1)+\gamma+q\delta$, $\nu=\gamma+\delta$ and $\beta+\gamma+\delta\equiv m\pmod{q-1}$, $m\in\{0,\ldots,q-2\}$. 
Let us write $$\bsb{f}=u^m\sum_{n\geq 0}a_n(t)v^n$$ with $a_n\in \FF_q[t,\theta]$ (this is possible after Proposition \ref{propositionEuexpansion} and 
the integrality of the coefficients of the $u$-expansions of $g,h$). Then, for all $n\geq 0$,
$$\deg_ta_n(t)\leq \nu\log_q\max\{1,n\}.$$

\item[(ii)] Let $\lambda$ be a positive real number. Let $\bsb{f}_1,\ldots,\bsb{f}_\sigma$ be a basis of monic monomials in $g,h,\bsb{E},\bsb{F}$ of
the $\mathcal{K}$-vector space $\MM^\dag_{\mu,\nu,m}(\mathcal{K})$. Let $x_1,\ldots,x_\sigma$ be polynomials of $\FF_q[t,\theta]$ with 
$\max_{0\leq i\leq \sigma}\deg_tx_i\leq \lambda$. Then, writing 
$$\bsb{f}=\sum_{i=1}^\sigma x_i\bsb{f}_i=u^m\sum_{n\geq 0}b_n(t)v^n$$ with $b_n\in \FF_q[t,\theta]$ with $0\leq m\leq q-2$, we have, for all $n\geq 0$:
$$\deg_tb_n\leq \lambda+\nu\log_q\max\{1,n\}.$$
\end{description}
\end{Proposition}

\noindent\emph{Proof.} 
Since by definition $\bsb{F}=\bsb{E}^{(1)}$, we have
$$\bsb{F}=u^q\sum_{n\geq 0}c_n^{(1)}v^{qn}=u\sum_{r\geq 0}d_rv^{r},$$ where $d_r=0$ if $q\nmid r-1$
and $d_r=c_{(r-1)/q}^{(1)}$ otherwise. Now, the operator $\tau$ leaves the degree in $t$ invariant. Therefore, by Proposition \ref{propositionEuexpansion}
$\deg_td_r\leq\log_q\max\{1,r/q\}\leq \log_q\max\{1,r\}$.

Let us consider the $u$-expansions:
$$\begin{array}{rclcrcl}g&=&\sum_{n\geq 0}\gamma_nv^n, & &\bsb{E}&=&u\sum_{n\geq 0}c_nv^n,\\
h&=&u\sum_{n\geq 0}\rho_nv^n,& & \bsb{F}&=&u\sum_{n\geq 0}d_nv^n,\end{array}$$ with $\gamma_n,\rho_n\in A$, $c_n,d_n\in\FF_q[t,\theta]$ for all $n$,
we can write:
$$\bsb{f}=u^{m'}\sum_{n\geq 0}\kappa_n v^n,$$
where $m'=\beta+\gamma+\delta$ and for all $n$, $\kappa_n=\sum\prod_x \gamma_{i_x}\prod_y \rho_{j_y}\prod_s c_{k_s}\prod_zd_{r_z}$, the sum being 
over the vectors of $\ZZ_{\geq 0}^{\alpha+\beta+\gamma+\delta}$ of the form
$$(i_1,\ldots,i_\alpha,j_1,\ldots,j_\beta,k_1,\ldots,k_\gamma,r_1,\ldots,r_\delta)$$ whose sum of entries is 
$n$, and with the four products running over
$x=0,\ldots,\alpha$, $y=0,\ldots,\beta$, $s=0,\ldots,\gamma$ and $z=0,\ldots,\delta$ respectively.
Since the coefficients of the $u$-expansions of $g,h$ do not depend on $t$ and $\gamma+\delta=\nu$, we obtain
$\deg_t\kappa_n\leq \nu\log_q\max\{1,n\}$.

If $m'=m+k(q-1)$ with $k\geq 0$ integer, and $0\leq m<q-1$. We can write 
$$\bsb{f}=u^{m'}\sum_{n\geq 0}c_n'v^n=u^{m}\sum_{n\geq 0}c_nv^n,$$ where $c_n=c_{n-k}'$, with the assumption that $c_{n-k}'=0$ if the index is negative.
The inequalities $\deg_tc_n'\leq\nu\log_q\max\{1,n\}$ for $n\geq 0$ imply that $\deg_tc_n$ is submitted to the same bound, proving the first part of the proposition.
The second part is a direct application of the first and ultrametric inequality.\CVD

\subsubsection{Dimensions of spaces}

\begin{Lemme}\label{usefulbunds}
We have, for all $m$ and $\mu,\nu\in\ZZ$ such that $\mu\geq(q+1)\nu\geq 0$,
$$\sigma(\mu,\nu)-\nu-1\leq \dim_{\mathcal{K}}\MM^\dag_{\mu,\nu,m}(\mathcal{K})\leq\sigma(\mu,\nu)+\nu+1,$$
where $$\sigma(\mu,\nu)=\frac{(\nu+1)\left(\mu-\frac{\nu(q+1)}{2}\right)}{q^2-1}.$$
Therefore, if $\mu>\frac{\nu(q+1)}{2}+q^2-1$, we have $\dim_{\mathcal{K}}\MM^\dag_{\mu,\nu,m}(\mathcal{K})>0$.
\end{Lemme}
\noindent\emph{Proof.} By \cite[p. 33]{Ge1}, we know that
$$\delta(k,m):=\dim_CM_{k,m}=\left\lfloor\frac{k}{q^2-1}\right\rfloor+\dim_CM_{k^*,m},$$
where $k^*$ is the remainder of the 
euclidean division of $k$ by $q^2-1$. In the same reference, it is also proved that
$\dim_CM_{k^*,m}=0$ unless $k^*\geq m(q+1)$, case where $\dim_CM_{k^*,m}=1$, so that, in all cases,
$0\leq\dim_CM_{k^*,m}\leq 1$. 

A basis of $\MM^\dag_{\mu,\nu,m}(\mathcal{K})$ is given by:
\begin{equation}\label{basis}(\bsb{b}_k)_{k=1,\ldots\dim\MM^\dag_{\mu,\nu,m}(\mathcal{K})}=(\phi_{i,s}\bsb{h}^s\bsb{E}^{\nu-s})_{s=0,\ldots,\nu,i=1,\ldots,\sigma(s)},\end{equation}
with, for all $s$, $(\phi_{i,s})_{i=1,\ldots,\sigma(s)}$ a basis of 
$M_{\mu-s(q-1)-\nu,m-\nu}$ (hence $\sigma(s)=\delta(\mu-s(q-1)-\nu,m-\nu)$). We have (taking into account the hypothesis on $\mu$ 
which implies $\mu-s(q-1)-\nu> 0$ for all $0\leq s\leq \nu$):
\begin{eqnarray*}
\dim\MM^\dag_{\mu,\nu,m}(\mathcal{K})&=&\sum_{s=0}^\nu\delta(\mu-\nu-s(q-1),m-\nu)\\
&=&\sum_{s=0}^\nu\left\lfloor\frac{\mu-s(q-1)-\nu}{q^2-1}\right\rfloor+\dim_CM_{(\mu-\nu-s(q-1))^*,m-\nu}.
\end{eqnarray*}
But 
$$\sum_{s=0}^\nu\frac{\mu-s(q-1)-\nu}{q^2-1}=\sigma(\mu,\nu).$$
Moreover, $\mu>\frac{\nu(q+1)}{2}+q^2-1$ if and only if $\sigma(\mu,\nu)>\nu+1$,
from which we deduce the lemma easily.\CVD

\subsubsection{Applying a variant of Siegel's Lemma}

We now prove the following:

\begin{Proposition}\label{functionsfmunum}
Let $\mu,\nu\in\ZZ_{\geq 0}$ be such that 
\begin{equation}\label{conditionmunu}
\mu\geq(q+1)\nu+2(q^2-1)
\end{equation} with $\nu\geq 1$, let $m$ be an integer in $\{0,\ldots,q-2\}$. There exists an integer $r>0$ such that 
\begin{equation}\label{boundforw}
r\leq 4q\mu\nu\log_q(\mu+\nu+q^2-1)+\nu
\end{equation} and, in $\widetilde{M}^{\leq \nu}_{r,m}$, 
a quasi-modular form $f_{\mu,\nu,m}$ such that 
\begin{equation}\label{boundfornu}
\frac{1}{q(q+1)}\mu\nu^2\log_q(\mu+\nu+q^2-1)\leq\nu_\infty(f_{\mu,\nu,m})\leq 4q\mu\nu^2\log_q(\mu+\nu+q^2-1).\end{equation}
\end{Proposition}

We will need the following variant of Siegel's Lemma whose proof can be found, for example, in \cite[Lemma 1]{LdM} (see also \cite{Bu}).
\begin{Lemme}\label{Siegel}
Let $U,V$ be positive integers, with $U<V$. Consider a system (\ref{system_siegel}) of 
$U$ equations with $V$ indeterminates:
\begin{equation}\label{system_siegel}
\sum_{i=1}^V a_{i,j}x_i=0,\quad (1\leq j\leq U)
\end{equation}
where the coefficients $a_{i,j}$ are elements of $K[t]$. Let $d$ be a non-negative integer such that
$\deg_ta_{i,j}\leq d$ for each $(i,j)$. Then, (\ref{system_siegel}) has a non-zero solution $(x_i)_{1\leq i\leq V}\in(K[t])^V$ with $\deg_tx_i\leq Ud/(V-U)$ for each $i=1,\ldots,V$.
\end{Lemme}
 
\noindent\emph{Proof of Proposition \ref{functionsfmunum}.} We apply Lemma \ref{Siegel} with the parameters $V=\dim\MM^\dag_{\mu,\nu,m}(\mathcal{K})$,
$U=\lfloor V/2\rfloor$.
We know that
$V>0$ because of (\ref{conditionmunu}) and Lemma \ref{usefulbunds}.

If $\bsb{f}=\bsb{b}_i$ as in (\ref{basis}), Writing
\begin{equation}\label{basisexpansion}
\bsb{b}_i=u^m\sum_{j\geq 0}a_{i,j}v^j,\quad a_{i,j}\in A[t]
\end{equation}
with $0\leq m< q-1$, Proposition \ref{propositionexpansions} says that for all $i$ and for all $j\geq 0$,
\begin{equation}\label{basisestimate}
\deg_ta_{i,j}\leq \nu\log_q\max\{1,j\}.
\end{equation}
Lemma \ref{Siegel} yields polynomials $x_1,\ldots,x_V\in K[t]$, not all zero, such that if we write 
\begin{equation}\label{deffunctionbsbf}
\bsb{f}=\sum_ix_i\bsb{b}_i=u^m\sum_{n\geq n_0}b_nv^n,\quad 0\leq m< q-1
\end{equation}
with $b_n\in K[t]$ for all $n$ and $b_{n_0}\neq 0$, we have the following properties.
The first property is the last inequality below:
\begin{eqnarray}
m+(q-1)n_0=\nu_\infty(\bsb{f})&\geq& m+(q-1)U\nonumber\\
&\geq &(q-1)(\sigma(\mu,\nu)-\nu-1)/2-1\nonumber\\
&\geq &\frac{(\nu+1)(\mu-\frac{\nu(q+1)}{2}-q^2+1)}{2(q+1)}-1\nonumber\\
&\geq&\frac{1}{4(q+1)}(\nu+1)\mu-1\label{minorationmultiplicity},
\end{eqnarray}
where we have applied Lemma \ref{usefulbunds} and (\ref{conditionmunu}).
The second property is that, in (\ref{deffunctionbsbf}),
\begin{equation}\label{boundbn}\deg_tb_n\leq2\nu(\log_q(\mu+\nu+q^2-1)+\log_q\max\{1,n\}),\quad n\geq 0,
\end{equation}
which follows from the following inequalities, with $d=\nu\log_q\max\{1,U\}$
\begin{eqnarray*}
\deg_tx_i&\leq & Ud/(V-U)\\
&\leq & \nu\log_q\max\{1,U\}\\
&\leq & \nu\log_q((\sigma(\mu,\nu)+\nu+1)/2)\\
&\leq & \nu(\log_q(\nu+1)+\log_q(\mu+q^2-1)-\log_q(q^2-1))\\
&\leq & 2\nu\log_q(\mu+\nu+q^2-1),
\end{eqnarray*}
and Proposition \ref{propositionexpansions}.

By Proposition \ref{Mdag}, 
we have $m+(q-1)n_0=\nu_\infty(\bsb{f})\leq \mu\nu$ so that $n_0\leq \frac{\mu\nu}{q-1}$,
where $n_0$ is defined in (\ref{deffunctionbsbf}).
Hence, by (\ref{boundbn}),
\begin{equation}\label{bounddegbn0}
\deg_tb_{n_0}\leq 4\nu\log_q(\mu+\nu+q^2-1).
\end{equation}

Lemma \ref{reduced} implies that for every integer $k$ such that  
\begin{equation}\label{boundonk}
k\geq \log_q(4\nu)+\log_q\log_q(\mu+\nu+q^2-1),
\end{equation}
the function $f_k:=\varepsilon(\bsb{f}^{(k)})$ satisfies $\nu_\infty(f_k)=\nu_\infty(\bsb{f}^{(k)})=q^k\nu_\infty(\bsb{f})$.
Let $k$ be satisfying (\ref{boundonk}). We have, by (\ref{minorationmultiplicity}), Proposition \ref{Mdag} and (\ref{ftauk}):
\begin{enumerate}
\item $f_k\in\widetilde{M}^{\leq \nu}_{\mu q^k+\nu,m}$,
\item $\left(\frac{(\nu+1)\mu q^k}{4(q+1)}-1\right)\leq\nu_\infty(f_k)\leq\mu\nu q^k$.
\end{enumerate}
Let us define the function 
$$\kappa(\mu,\nu):=\lfloor\log_q(4\nu)+\log_q\log_q(\mu+\nu+q^2-1)\rfloor+1$$
and write:
$f_{\mu,\nu,m}:=f_{\kappa(\mu,\nu)}$. This Drinfeld quasi-modular form satisfies the properties announced in the proposition.\CVD

\subsection{Proof of Theorem \ref{secondtheorem}}

Let $f$ be a Drinfeld quasi-modular form of weight $w$ and depth $l$. We can assume, without loss of generality, that $f$, as a polynomial in $E,g,h$ with coefficients in $C$,
it is an irreducible polynomial. We can also assume, by Gekeler, \cite[Formula (5.14)]{Ge} and, \cite[Theorem 1.4]{BP2},
that $l>q$.

Let $W$ be a real number $\geq 1$ and 
let $\alpha$ be the function of a real variable defined, for $\mu\geq 0$, by $\alpha(\mu)=\mu l\log_q(\mu+Wl+q^2-1)$; 
we have $\alpha(\mu+1)\leq 2\alpha(\mu)$.
Since (the dash $'$ is the derivative) $\alpha'(\mu)\geq l\log_q(Wl+q^2-1)>1$
for all $l\geq q$ and $\mu\geq0$, for all $w\geq 0$ integer, there exist $\mu\in\ZZ_{\geq 0}$ such that 
\begin{equation}\label{ulbound}
\alpha(\mu)\leq w<\alpha(\mu+1),
\end{equation}
and we choose one of them, for example the biggest one. 
Let us suppose that (\ref{condition}) holds and, at once, set
$$\nu=Wl,$$
with 
$$W=q(2+4(q+1))=2 q (3 + 2 q).$$

We define $\beta(l)$ to be the right hand side of (\ref{condition}), as a function of $l\geq q$.
Condition (\ref{condition}) implies $$\mu\geq\frac{\beta(l)}{2l\log_q(\mu+Wl+q^2-1)}.$$
Since $\log_q(x)\leq 2x^{1/2}$ for all $x\geq 1$ and $q\geq2$, we get
$$(\mu+Wl+q^2-1)^{3/2}\geq\frac{\beta(l)}{4l},$$
that is,
$$\mu\geq\left(\frac{\beta(l)}{4l}\right)^{2/3}-Wl-q^2+1.$$ But replacing $\beta(l)$ by its value yields $\mu\geq(q+1)\nu+2(q^2-1)$, which is the condition (\ref{conditionmunu}) needed to apply Proposition \ref{functionsfmunum}.

Let us write $\mathcal{L}:=\log_q(\mu+\nu+q^2-1)$ so that $\alpha(\mu)=\mu l\mathcal{L}$. By Proposition \ref{functionsfmunum}, there exists a 
form $f_{\mu,\nu,m}\in\widetilde{M}^{\leq \nu}_{r,m}$ such that $l(f_{\mu,\nu,m})\leq\nu$ and
\begin{equation}\label{pmatrix}
\begin{array}{rcccl}
& &w(f_{\mu,\nu,m})&\leq&4(q+1)\mu\nu\mathcal{L}\\
(q(q+1))^{-1}\mu\nu^2\mathcal{L}&\leq &\nu_\infty(f_{\mu,\nu,m})&\leq&4q\mu\nu^2\mathcal{L}
\end{array}
\end{equation}

We have two cases.

\medskip

\noindent\emph{Case (i)}. If $f|f_{\mu,\nu,m}$, then 
\begin{equation}\label{case1}
\nu_\infty(f)\leq\nu_\infty(f_{\mu,\nu,m})\leq4q\mu\nu^2\mathcal{L}.
\end{equation}

\medskip

\noindent\emph{Case (ii)}. If $f\nmid f_{\mu,\nu,m}$, then $\rho:=\mathbf{Res}_E(f,f_{\mu,\nu,m})$
is a non-zero modular form, whose weight $w(\rho)$ and type $m(\rho)$ can be computed with the help of \cite[Lemma 2.5]{BP2}
(we do not need an explicit computation of $m(\rho)$):
\begin{eqnarray}
w(\rho)&=&w\nu+w(f_{\mu,\nu,m})l-2l\nu\nonumber\\
&\leq & w\nu+4l(q+1)\mu\nu\mathcal{L}-2l\nu\nonumber\\
&\leq &\nu(w+4(q+1)\mu l\mathcal{L})\nonumber\\
&< &\nu(\alpha(\mu+1)+4(q+1)\mu l\mathcal{L})\nonumber\\
&<&\nu(2\alpha(\mu)+4(q+1)\mu l\mathcal{L})\nonumber\\
&<&(2+4(q+1))\nu\mu l\mathcal{L}.\label{boundwrho}
\end{eqnarray}

Let us suppose that $\nu_\infty(f)>(q(q+1))^{-1}\mu\nu^2\mathcal{L}$. Then, by B\'ezout identity for the resultant and (\ref{pmatrix}), $\nu_\infty(\rho)\geq(q(q+1))^{-1}\mu\nu^2\mathcal{L}$. At the same time, by Gekeler, \cite[Formula (5.14)]{Ge}, $\nu_\infty(\rho)\leq\frac{w(\rho)}{q+1}$, yielding the inequality $W<q(2+4(q+1))$ which 
is contradictory with the definition of $W$.

Therefore, in case (ii), we have that $\nu_\infty(f)\leq4q\mu\nu^2\mathcal{L}$. Ultimately, we have shown that, in both
cases  (i), (ii),
\begin{eqnarray*}
\nu_\infty(f)&\leq&4q\mu\nu^2\mathcal{L}\\
&\leq&4q\mu W^2l^2\mathcal{L}\\
&\leq &4qW^2lw,
\end{eqnarray*}
which is the estimate (\ref{finalupperbound}).\CVD

\begin{Remarque} {\em The dependence on $l$ in condition (\ref{condition}) can be 
relaxed, adding conditions on $q$. For all $\epsilon>0$ there exists a constant $c>0$ such that for all $q>c$, assuming that $w\gg_{\epsilon}l^{2+\epsilon}$, then, the inequality (\ref{finalupperbound}) holds. We do not report the proof of this fact here.}\end{Remarque}

\medskip

\noindent\emph{Acknowledgement.} The author is indebted with V. Bosser for several discussions on these topics and a careful reading of the first versions of this paper.


\begin{thebibliography}{99}

\bibitem{An} G. Anderson. {\em $t$-motives}, Duke Math. J. 53 (1986), 457-502. 

\bibitem{ABP} G. Anderson, D. Brownawell \& M. Papanikolas, {\em Determination of 
the algebraic relations among special $\Gamma$-values in positive characteristic,} Ann. of Math. 160 
(2004), 237-313. 



\bibitem{BP} V. Bosser \& F. Pellarin. {\em Differential properties of Drinfeld quasi-modular forms.}
Int. Math. Res. Notices. Vol. 2008.

\bibitem{BP2} V. Bosser \& F. Pellarin. {\em On certain families of Drinfeld quasi-modular forms.}
to appear in J. Number Theory  (2009), {\tt doi:10.1016/j.jnt.2009.04.014}.


\bibitem{BM} D. Brownawell \& D. Masser. {\em Multiplicity estimates for analytic functions I.} J. Reine angew. Math. 314 (1980), pp. 200-216.

\bibitem{BM2} D. Brownawell \& D. Masser. {\em Multiplicity estimates for analytic functions II.} Duke Math. J. Volume 47, Number 2 (1980), pp. 273-295.

\bibitem{Bu} P. Bundschuh. {\em Transzendenzma\ss e in K\"orpern, Laurentreihen.} J. reine angew. Math. 299-300 (1978), pp. 411-432.

\bibitem{ChaPa} Chieh-Yu Chang \& M. Papanikolas. {\em Algebraic relations among periods and logarithms of rank $2$ Drinfeld modules}. Preprint.\\ {\tt http://arxiv.org/abs/0807.3157}



\bibitem{FP} J. Fresnel, \& M. van der Put. {\em Rigid Analytic Geometry and its Applications.} Birkh\"auser, 
Boston (2004).




\bibitem{gekeler:compositio} E.-U. Gekeler. {\em Quasi-periodic functions and Drinfeld modular forms.}
Compositio Math. t. 69 No. 3 p. 277-293 (1989).

\bibitem{Ge}  E.-U. Gekeler. {\em
On the coefficients of Drinfeld modular forms.}
Invent. Math. 93, No.3, 667-700 (1988).

\bibitem{Ge1} E.-U. Gekeler.  
{\em Lectures on Drinfeld Modular Forms.}
(Written by I. Longhi) CICMA Lecture Notes No. 4 (1999).


\bibitem{Go2} D. Goss. {\em Basic structures of function field arithmetic.} Ergebnisse der Mathematik und ihrer Grenzgebiete, 35. Springer-Verlag, Berlin, 1996.


\bibitem{KK1} M. Kaneko \& M. Koike. {\em On Modular Forms Arising from a Differential 
Equation of Hypergeometric Type.} Ramanujan Journ. 7, pp. 145-164, (2003).



\bibitem{KZ1} M. Kaneko \& D. Zagier, {\em Supersingular $j$-invariants, Hypergeometric series, and Atkin's orthogonal polynomials,} AMS/IP Studies in Advanced Mathematics 7, 97-126, (1998). 

\bibitem{LdM} A. Lasjaunias \& B. de Mathan. {\em Thue's Theorem in positive characteristic.} J. reine angew. Math. 473, 195-206, (1996).




\bibitem{Nes}
Yu. V. Nesterenko. {\em Estimate of the orders of the zeroes of functions of a certain class, and their application in the theory of transcendental numbers.} Izv. Akad. Nauk SSSR Ser. Mat. 41 (1977), no. 2, pp. 253-284, 477, Math. USSR Izv. 11 (1977), pp. 239-270.

\bibitem{NP} Yu. V. Nesterenko et P. Philippon \'editeurs. 
{\em Introduction to algebraic independence theory.} Lecture Notes in Mathematics 1752, Springer (2001).

\bibitem{Pa} M. A. Papanikolas. {\em Tannakian duality for Anderson-Drinfeld
motives and algebraic independence of Carlitz logarithms},
Invent. Math. 171, 123-174 (2008).

\bibitem{PF} F. Pellarin.
{\em 
La structure diff\'erentielle de l'anneau des formes quasi-modulaires pour ${\bf SL}_2({\bf Z})$.} Journal de th\'eorie des nombres de Bordeaux, 18 no. 1 (2006), p. 241-264 

\bibitem{Bourbaki} F. Pellarin. {\em Aspects de l'ind\'ependance alg\'ebrique en caract\'eristique non nulle.}
S\'eminaire Bourbaki Mars 2007 
59me ann\'ee, 2006-2007, no. 973.



\bibitem{Stiller} P. Stiller. {\em Special values of Dirichlet series, monodromy, and the periods of automorphic
forms.} Mem. Am. Math. Soc. 49 (299), (1984).

\end{thebibliography}
\end{document}